\documentclass[11pt, a4paper]{article}
\usepackage{CJKutf8}
\usepackage{geometry}
\usepackage{amssymb}
\usepackage{amsmath}
\usepackage{amsthm}
\usepackage{amsfonts}
\usepackage{mathtools}
\usepackage{mathrsfs}
\usepackage{bm} 

\usepackage{graphicx}
\usepackage{dcolumn}
\usepackage[ruled,lined,algonl]{algorithm2e}
\usepackage{color}
\usepackage{autobreak}
\allowdisplaybreaks

\theoremstyle{plain}
\newtheorem{theorem}{Theorem}
\newtheorem{corollary}{Corollary}

\newtheorem{proposition}{Proposition}

\theoremstyle{definition}
\newtheorem{definition}{Definition}

\theoremstyle{remark}

\newcommand{\pset}[1]{\mathcal{#1}}

\DeclareMathOperator{\numer}{Numer}
\DeclareMathOperator{\denom}{Denom}

\newcommand{\uvar}{\bm{u}}
\newcommand{\xvar}{\bm{x}}

\usepackage{authblk}

\begin{document}
\begin{CJK}{UTF8}{gbsn}

\title{Analysis of stability and bifurcation for two heterogeneous triopoly games with the isoelastic demand}

\author{Xiaoliang Li\thanks{Corresponding author: xiaoliangbuaa@gmail.com}}
\affil{School of Finance and Trade, Dongguan City College, Dongguan, P. R. China}

%
%
%

\maketitle

\begin{abstract}

In this paper, we investigate two heterogeneous triopoly games where the demand function of the market is isoelastic. The local stability and the bifurcation of these games are systematically analyzed using the symbolic approach proposed by the author. The novelty of the present work is twofold. On one hand, the results of this paper are analytical, which are different from the existing results in the literature based on observations through numerical simulations.  In particular, we rigorously prove the existence of double routes to chaos through the period-doubling bifurcation and through the Neimark-Sacker bifurcation. On the other hand, for the special case of the involved firms having identical marginal costs, we acquire the necessary and sufficient conditions of the local stability for both models. By further analyzing these conditions,  it seems that that the presence of the local monopolistic approximation (LMA) mechanism might have a stabilizing effect for heterogeneous triopoly games with the isoelastic demand.\\

\par\textbf{Keywords:} triopoly games; heterogeneous firms; dynamics; stability; bifurcations; symbolic computation
\end{abstract}

\section{Introduction}

Different from a competitive market that possesses a large number of relatively small companies producing homogeneous products and competing with each other, an oligopoly is a market supplied by only a few firms.  It is well known that Cournot developed the first formal theory of oligopoly in \cite{Cournot1838R}, where players are supposed to have the naive expectations that their rivals produce the same quantity of output as in the immediately previous period. Since a firm's decision as well as the reactions of its competitors can change the price and further influence their profits, an oligopolistic market is much more complex than a competitive market although only a few firms are involved. 

Puu \cite{Puu1991C} followed Cournot's seminal work by considering the isoelastic demand function and setting the price simply the reciprocal of the total supply of the market. He introduced a game where the players also have the naive expectations, and then showed that the evolution of its dynamics can develop a period-doubling bifurcation and finally to chaos. It was pointed out in \cite{Hommes2013B} that the expectations and the behavioral rationalities play ambitious roles in modeling economic phenomena. As far as we know, there mainly exist four types of expectations in the current literature: the naive \cite{Puu1991C}, the boundedly rational \cite{Elsadany2012C}, the adaptive \cite{Bischi2001E} and the local monopolistic approximation (or LMA for short) \cite{Tuinstra2004A}.

The monopoly, a market ruled by a unique firm, is the simplest oligopoly, but may exhibit quite complex dynamic behaviors such as quasi-period and chaos \cite{Puu1995T, Naimzada2008C, Matsumoto2021D, AlHdaibat2015O}. However, in the real world, it is more possible that a market supplied by two, three, or even four companies, which is called a duopoly, a triopoly, or a quadropoly, respectively. Among them, the duopoly is most intensively studied, see, e.g., 
\cite{Agiza1999O, Ahmed2000O, Askar2016T, Canovas2018O,Govaerts2008S,Kopel1996S, Li2020N,Ma2013T,Peng2016A,Sarafopoulos2015Ca,Shi2015C,Wu2010C}. Furthermore, the reader can refer to \cite{Andaluz2017N,Elsadany2012C,Ma2011C,Ma2013C,Tramontana2012H,Tu2017C} for the explorations of triopoly games, which are also the topic of our concern in this paper. Moreover, there are relatively few studies on quadropoly models, and the reader may refer to \cite{Elsadany2013C,Matouk2017N}.

In the present work, we investigate two distinct heterogeneous triopoly games with the isoelastic demand function. One game was first proposed in \cite{Tramontana2012H}, where an adoptive firm, a naive firm and a boundedly rational firm compete with each other. By numerical simulations, it was observed that this game could exhibit both the period-doubling bifurcation and the Neimark-Sacker bifurcation as the parameters vary. The other game \cite{Naimzada2015T} is similar but with the first adoptive player replaced by an LMA player. For this game, it was also discovered that there exist the same two possible routes to chaos, through the period-doubling bifurcation and through the Neimark-Sacker bifurcation.

We try to apply the symbolic approach proposed by the author of this paper and his coworker in \cite{Li2014C}, to algebraically analyze the stability and the bifurcation of the two models mentioned above. Consider the parametric system of the following form.
\begin{equation}\label{eq:par-sys}
\left\{
\begin{array}{l}\smallskip
F_1(\uvar,\xvar)=0,\ldots,F_n(\uvar,\xvar)=0,\\ \smallskip
N_1(\uvar,\xvar)\neq 0,\ldots,N_m(\uvar,\xvar)\neq 0,\\ \smallskip
P_1(\uvar,\xvar)> 0,\ldots,P_s(\uvar,\xvar)> 0,\\
P_{s+1}(\uvar,\xvar)\geq 0,\ldots,P_{s+t}(\uvar,\xvar)\geq 0,
\end{array}
\right.
\end{equation}
where $\uvar=(u_1,u_2,\ldots,u_d)$ are the parameters, $\xvar=(x_1,x_2,\ldots,x_n)$ are the variables, and $F_1,\ldots,F_n$, $N_1,\ldots,N_m$, $P_1,\ldots,P_{s+l}$ are all polynomials. The symbolic approach in \cite{Li2014C} permits us to obtain the necessary and sufficient condition that system \eqref{eq:par-sys} has a given number of distinct real solutions. The basic idea of the approach is to first transform system \eqref{eq:par-sys} into a new one in a single variable and then analyze the real solutions of the resulting system. This approach is different from those based on numerical computations. It can be used to establish exact and rigorous results and thus are more adequate for the theoretical study of economic models.

Take the following system as an illustrative example.
\begin{equation}\label{ex:with-par}
\left\{
\begin{array}{l}\smallskip
x^3-uy^2=0,\\ \smallskip y^2-2\,x-1=0, \\ \smallskip
x-y\neq 0,\\
y+s >0,
\end{array}
\right.
\end{equation}
where $u,s$ are parameters. The first step is transforming the system into a  univariate system as follows such that the original and the resulting have the same number of distinct real solutions.
\begin{equation}\label{eq:uni-sys}
\left\{\begin{split}
&x^6-(4\,u+3)x^4-18\,ux^3+(4\,u^2-26\,u+3)x^2\\
        &+(4\,u^2-14\,u)x+u^2-2\,u-1=0,\\
&(-J+Is)I>0,
\end{split}\right.
\end{equation}
where 
\begin{equation*}
\begin{split}
    &I=-3\,x^2-8\,x+2\,u-5,\\
    &J=-x^3-6\,x^2+(2\,u-7)x+u-2.	
\end{split}	
\end{equation*}
We can then decompose the parameter space $(u,s)$ into regions such that on each of them the number of distinct real solutions of \eqref{eq:uni-sys}, or \eqref{ex:with-par} equivalently, is invariant. The computational results finally conclude: provided that
$$u(32\,u-27)(u^2-2\,u-1)R\neq 0,$$
where $R=s^6-3\,s^4-8\,us^2+3\,s^2-1$, the number of distinct real
solutions of system \eqref{ex:with-par} should be
\begin{itemize}
  \item $0$ if and only if $R<0$ and $s<0$;
  \item $1$ if and only if $R>0$;
  \item $2$ if and only if $R<0$ and $s>0$.
\end{itemize}
The reader can refer to \cite{Li2014C} for more details.

The novelty of this paper is twofold. For the two triopoly games considered in this paper, the existing investigations are mainly built on observations through numerical simulations. By contrast, our computations are exact and rigorous, which not only confirm the existing results but also provide a solid foundation for the theoretical development of triopoly games.  Moreover, for the special case of the involved firms having identical marginal costs, we acquire the necessary and sufficient conditions of the local stability for both models. By further analyzing these conditions, we show that the stable regions of the two models would be the same if the proportion parameter takes the value $2/3$. Our analysis confirms the statement in \cite{Bischi2007O, Naimzada2009C, Cavalli2015N} that the presence of the LMA mechanism might have a stabilizing effect and suggests that one may extend a similar conclusion to the case of generic size heterogeneous oligopoly models with the isoelastic demand.

The rest of this paper is structured as follows. In Section \ref{sec:models}, we revisit the construction of the two models. In Section \ref{sec:stability}, the local stability is thoroughly investigated, and the stable regions of the two games are compared analytically for the special case of identical marginal costs. In Section \ref{sec:pd-bifur} and Section \ref{sec:ns-bifur}, by virtue of the algebraic criteria proposed by Wen and others \cite{Wen2005C,Wen2008A}, we obtain the complete conditions of the occurrence of the period-doubling bifurcation and the Neimark-Sacker bifurcation.  The paper is concluded with some remarks in Section \ref{sec:con}.

\section{Models}\label{sec:models}

Let us consider a market served by three firms producing homogeneous products. We use $x(t)$, $y(t)$ and $z(t)$ to denote the outputs at period $t$ of the three firms, respectively. Assume that the demand function is isoelastic, which is founded on the hypothesis that the consumers have the Cobb-Douglas utility function. Hence, the price of the product should be
$$p(Q)=\frac{1}{Q}=\frac{1}{x+y+z},$$
where $Q=x+y+z$ is the total supply. Moreover, the marginal cost of each firm is supposed to be constant. To be exact, the cost functions of the three firms are set to be $c_1x$, $c_2y$ and $c_3z$ ($c_1,c_2,c_3>0$), respectively.  

Under the above assumptions, the first firm should have the profit function
$$\Pi_1(x,y,z)=\frac{x}{x+y+z}-c_1x.$$
Evidently, its marginal profit would be
$$\phi_1(x,y,z)=\frac{y+z}{(x+y+z)^2}-c_1.$$
The profits $\Pi_2,\Pi_3$ and the marginal profits $\phi_2,\phi_3$ of the second and the third firm can be similarly obtained. 

Suppose that the first firm is an \emph{adaptive} player. At each period $t+1$, the first firm decide the output $x(t+1)$ according to the previous output $x(t)$ and its expectations of the other two competitors. It is supposed that firm $1$ naively expects that at period $t+1$ firm $2$ and $3$ would produce the same quantity as at period $t$. Therefore, the first firm could calculate the optimal output $x_{opt}$ to maximize its expected profit $\Pi_1(x_{opt},y(t),z(t))$, i.e., such that
$$\phi_1(x_{opt},y(t),z(t))=0.$$
That is
$$x_{opt}=\sqrt{\frac{y(t)+z(t)}{c_1}}-y(t)-z(t).$$
The adaptive decision mechanism for firm $1$ is that it choose the output $x(t+1)$ proportionally to be
$$x(t+1)=(1-l)x(t)+lx_{opt}=(1-l)x(t)+l\left[\sqrt{\frac{y(t)+z(t)}{c_1}}-y(t)-z(t)\right],$$
where $l\in[0,1]$ is a parameter controlling the proportion. 

The second firm is simply a \emph{naive} player, who expects that the productions of the rivals would be equal to those of the last period. Thus, in order to gain the maximal profit, it would choose the best response to produce
$$y(t+1)=\sqrt{\frac{x(t)+z(t)}{c_2}} - x(t) - z(t).$$

Our third firm is assumed to be a \emph{boundedly rational} player, who adopts the so-called gradient adjustment mechanism, i.e.,
$$z(t+1)=z(t)+kz(t)\phi_3(x(t),y(t),z(t))=z(t)+kz(t)\left[\frac{x(t)+y(t)}{\big(x(t)+y(t)+z(t)\big)^2}-c_3\right],$$
where $k>0$ is a parameter controlling the adjustment speed. In other words, the third firm increases/decreases its output according to the information given by the marginal profit of the last period. It is worth noting that the adjustment speed depends upon not only the parameter $k$ but also the size of the firm $z(t)$.

In short,  the dynamics of the above triopoly game could be described by the following $3$-dimensional discrete dynamic system, which is denoted by $T_{ANB}(x,y,z)$.
\begin{equation}\label{eq:iter-abg}
	T_{ANB}(x,y,z): 
	\left\{\begin{split}
	&x(t+1)=(1-l)x(t)+l\left[\sqrt{\frac{y(t)+z(t)}{c_1}}-y(t)-z(t)\right],\\
	&y(t+1)=\sqrt{\frac{x(t)+z(t)}{c_2}} - x(t) - z(t),\\
	&z(t+1)=z(t)+kz(t)\left[\frac{x(t)+y(t)}{\big(x(t)+y(t)+z(t)\big)^2}-c_3\right].
	\end{split}
	\right.
\end{equation}

The above game was first introduced by Tramontana and Elsadany in \cite{Tramontana2012H}. Motivated by their work, a distinct heterogeneous triopoly game was investigated in \cite{Naimzada2015T}, where the first firm is replaced with a firm that has further reduced rationality. This firm adopts the so-called \emph{local monopolistic approximation} (LMA) mechanism \cite{Tuinstra2004A}, which is also a boundedly rational adjustment process. In this process, the player just has limited knowledge of the demand function. In particular, the firm can observe the current market price $p(t)$ and the corresponding total supply $Q(t)$ and is able to correctly estimate the slope $p'(Q(t))$ of the price function around the point $(p(t),Q(t))$. Then, the firm uses such information to conjecture the demand function and expect the price at period $t+1$ to be
$$p^e(t+1)=p(Q(t))+p'(Q(t))(Q^e(t+1)-Q(t)),$$
where $Q^e(t+1)$ represents the expected aggregate production at period $t+1$. Moreover, firm $1$ is also assumed to use the naive expectations of its rivals, i.e., $Q^e(t+1)=x(t+1)+y(t)+z(t)$. Thus, we have that 
$$p^e(t+1)=\frac{1}{Q(t)}-\frac{1}{Q^2(t)}(x(t+1)-x(t)).$$
The expected profit of the first firm is
$$\Pi^e_1(t+1)=p^e(t+1)x(t+1)-c_1x(t+1).$$
To maximize the expected profit, firm $1$ chooses its output at period $t+1$ to be
$$x(t+1)=\frac{2\,x(t)+y(t)+z(t)-c_1\big(x(t)+y(t)+z(t)\big)^2}{2}.$$
Therefore, the new model can be described by the following $3$-dimensional discrete dynamic system, which is denoted by $T_{LNB}(x,y,z)$.

\begin{equation}\label{eq:iter-lbg}
	T_{LNB}(x,y,z): 
	\left\{\begin{split}
	&x(t+1)=\frac{2\,x(t)+y(t)+z(t)-c_1\big(x(t)+y(t)+z(t)\big)^2}{2},\\
	&y(t+1)=\sqrt{\frac{x(t)+z(t)}{c_2}} - x(t) - z(t),\\
	&z(t+1)=z(t)+kz(t)\left[\frac{x(t)+y(t)}{\big(x(t)+y(t)+z(t)\big)^2}-c_3\right].
	\end{split}
	\right.
\end{equation}

\section{Local Stability}\label{sec:stability}

In order to determine the equilibria of the two models introduced in Section \ref{sec:models}, we set $x(t+1)=x(t)=x$, $y(t+1)=y(t)=y$ and $z(t+1)=z(t)=z$. Then we have that, for the $T_{ANB}$ model, 
\begin{equation}
	\left\{\begin{split}
	&\sqrt{\frac{y+z}{c_1}}-(x+y+z)=0,\\
	&\sqrt{\frac{x+z}{c_2}}-(x+y+z)=0,\\
	&kz\left[-c_3+\frac{x+y}{(x+y+z)^2}\right]=0,
	\end{split}
	\right.
\end{equation}
and for the $T_{LNB}$ model,
\begin{equation}
	\left\{\begin{split}
	&y+z-c_1(x+y+z)^2=0,\\
	&\sqrt{\frac{x+z}{c_2}}-(x+y+z)=0,\\
	&kz\left[-c_3+\frac{x+y}{(x+y+z)^2}\right]=0.
	\end{split}
	\right.
\end{equation}
It is obvious that both of them can be solved by two equilibria:
$$E_1=\left(\frac{c_2}{(c_1+c_2)^2},\frac{c_1}{(c_1+c_2)^2},0\right),$$
$$E_2=\left(\frac{2(c_2+c_3-c_1)}{(c_1+c_2+c_3)^2}, \frac{2(c_1+c_3-c_2)}{(c_1+c_2+c_3)^2},\frac{2(c_1+c_2-c_3)}{(c_1+c_2+c_3)^2}\right).$$
It should be noted that the origin $(0,0,0)$ is not an equilibrium since it is not defined for the iteration maps. From the economic point of view, only the equilibrium $E_2$ is of our concern because the state $E_1$, where the third firm is out of the market forever, is unstable when $E_2$ is strictly positive as mentioned in \cite{Tramontana2012H}. 

The Jacobian matrix takes the form
\begin{equation*}
J = \left[\begin{matrix}
	\frac{\partial x(t+1)}{\partial x(t)} & \frac{\partial x(t+1)}{\partial y(t)} & \frac{\partial x(t+1)}{\partial z(t)}\\
	\frac{\partial y(t+1)}{\partial x(t)} & \frac{\partial y(t+1)}{\partial y(t)} & \frac{\partial y(t+1)}{\partial z(t)}\\
	\frac{\partial z(t+1)}{\partial x(t)} & \frac{\partial z(t+1)}{\partial y(t)} & \frac{\partial z(t+1)}{\partial z(t)}\\
\end{matrix}\right].
\end{equation*}
We use $J_{ABG}$ and $J_{LBG}$ to denote the Jacobian matrices of our two models, respectively. Then,
\begin{equation*}
J_{ABG} = \left[\begin{matrix}
1-l & l\left[\frac{1}{2c_1}\left(\frac{y+z}{c_1}\right)^{-\frac{1}{2}}-1 \right] & l\left[\frac{1}{2c_1}\left(\frac{y+z}{c_1}\right)^{-\frac{1}{2}}-1 \right] \\
\frac{1}{2c_2}\left(\frac{x+z}{c_2}\right)^{-\frac{1}{2}}-1 & 0 & \frac{1}{2c_2}\left(\frac{x+z}{c_2}\right)^{-\frac{1}{2}}-1 \\
\frac{kz(z-x-y)}{(x+y+z)^3} & \frac{kz(z-x-y)}{(x+y+z)^3} & 1-kc_3+\frac{k(x+y)(x+y-z)}{(x+y+z)^3}\\
\end{matrix}\right]
\end{equation*}
and
\begin{equation*}
J_{LBG} = \left[\begin{matrix}
1-c_1(x+y+z) & 1/2-c_1(x+y+z) & 1/2-c_1(x+y+z) \\
\frac{1}{2c_2}\left(\frac{x+z}{c_2}\right)^{-\frac{1}{2}}-1 & 0 & \frac{1}{2c_2}\left(\frac{x+z}{c_2}\right)^{-\frac{1}{2}}-1 \\
\frac{kz(z-x-y)}{(x+y+z)^3} & \frac{kz(z-x-y)}{(x+y+z)^3} & 1-kc_3+\frac{k(x+y)(x+y-z)}{(x+y+z)^3}\\
\end{matrix}\right].
\end{equation*}

For continuous dynamic systems, one may use the first method of Lyapunov with the technique of linearization. A similar procedure can be used to linearize the discrete system around the equilibrium and to determine the stability conditions via the eigenvalues of the Jacobian matrix of the system. The following proposition permits us to identify the stability of an equilibrium.

\begin{proposition}
Let $E$ be an equilibrium of a discrete dynamic system. If all the eigenvalues $\lambda_i$ of the Jacobian matrix at $E$ lie in the open unit disk, i.e., $|\lambda_i|<1$ for all $i$, then $E$ is locally stable. Moreover, if the Jacobian matrix at $E$ has at least one eigenvalue $\lambda_0$ outside the open unit disk, i.e., $|\lambda_0|>1$, then $E$ is unstable.
\end{proposition}

Let $A$ be the characteristic polynomial of the Jacobian matrix  $J$. The eigenvalues of $J$ are simply the roots of the polynomial $A$ for $\lambda$. So the problem of stability analysis can be reduced to that of determining whether all the roots of $A$ lie in the open unit disk $|\lambda|<1$.  To the best of our knowledge, in addition to the Routh-Hurwitz criterion \cite{Oldenbourg1948T} generalized from the corresponding criterion for continuous systems, there are two other criteria, the Schur-Cohn criterion \cite[pp.\,246--248]{Elaydi2005U} and the Jury criterion \cite{Jury1976I}, available for discrete dynamical systems. In what follows, we provide a short review of the Schur-Cohn criterion.

\begin{proposition}[Schur-Cohn Criterion]
For a $n$-dimensional discrete dynamic system, assume that the characteristic polynomial of its Jacobian matrix is
\begin{equation*}
A= \lambda^n + a_{n-1}\lambda^{n - 1} + \cdots + a_0.
\end{equation*}
Consider the
sequence of determinants $D^\pm_1$, $D^\pm_2$, $\ldots$, $D^\pm_n$,
where
\begin{equation*}
\begin{split}
  D^{\pm}_i =&\left| \left(
\begin{array}{ccccc}
1&a_{n-1}&a_{n-2}&\cdots&a_{n-i+1}\\
0&1&a_{n-1}&\cdots&a_{n-i+2}\\
0&0&1&\cdots&a_{n-i+3}\\
\vdots&\vdots&\vdots&\ddots&\vdots\\
0&0&0&\cdots&1\\
\end{array}
\right)\pm\left(
\begin{array}{ccccc}
a_{i-1}&a_{i-2}&\cdots&a_{1}&a_0\\
a_{i-2}&a_{i-3}&\cdots&a_{0}&0\\
\vdots&\vdots&\ddots&\vdots&\vdots\\
a_{1}&a_0&\cdots&0&0\\
a_0&0&\cdots&0&0\\
\end{array}
\right)\right|.
\end{split}
\end{equation*}
The characteristic polynomial $A$ has all its roots inside the unit
open disk if and only if
\smallskip
\begin{enumerate}
\item $A(1)>0$ and $(-1)^nA(-1)>0$,

\item  $D^\pm_1>0, D^\pm_3>0, \ldots, D^\pm_{n-3}>0, D^\pm_{n-1}>0$
(when $n$ is even), or\\[2pt]
\smallskip $D^\pm_2>0, D^\pm_4>0, \ldots, D^\pm_{n-3}>0,
D^\pm_{n-1}>0$ (when $n$ is odd).
\end{enumerate}
\end{proposition}

\begin{corollary}
	Consider a $3$-dimensional discrete dynamic system with the characteristic polynomial of its Jacobian matrix of the form 
	$$A=\lambda^3+a_2\lambda^2+a_1\lambda+a_0.$$
	An equilibrium $E$ is locally stable if and only if the following conditions are satisfied at $E$.
	\begin{equation}
		\left\{\begin{split}
			&1+a_2+a_1+a_0>0,\\
			&1-a_2+a_1-a_0>0,\\
			&-a_0^2-a_0a_2+a_1+1>0,\\
			&-a_0^2+a_0a_2-a_1+1>0.
		\end{split}\right.
	\end{equation}
\end{corollary}

According to the above corollary, the problem of determining the stability of the $T_{ANB}$ model is transformed into determining the existence of real solutions for the following system.

\begin{equation}\label{eq:radical-abg}
	\left\{\begin{split}
	&\sqrt{\frac{y+z}{c_1}}-(x+y+z)=0,\\
	&\sqrt{\frac{x+z}{c_2}}-(x+y+z)=0,\\
	&kz\left[-c_3+\frac{x+y}{(x+y+z)^2}\right]=0,\\
	&1+a_2+a_1+a_0>0,\\
	&1-a_2+a_1-a_0>0,\\
	&-a_0^2-a_0a_2+a_1+1 >0,\\
	&-a_0^2+a_0a_2-a_1+1>0,\\
	&x>0,~y>0,~z>0,\\
	&k>0,~1\geq l\geq 1,~c_1>0,~c_2>0,~c_3>0,
	\end{split}
	\right.
\end{equation}
where $x>0$, $y>0$, $z>0$ are added to limit our focus on $E_2$, and $k>0$, $1\geq l\geq 1$, $c_1>0$, $c_2>0$, $c_3>0$ are generated from the economic meanings of all the involved parameters. It is worth noticing that $a_0$, $a_1$ and $a_2$ are extremely complex rational functions, which involve the radical expressions $\sqrt{\frac{y+z}{c_1}}$ and $\sqrt{\frac{x+z}{c_2}}$. 

We herein try to determine the existence of real solutions for system \eqref{eq:radical-abg} by using the method proposed in \cite{Li2014C}, which is suitable only for polynomial systems. The obvious obstacle for us is that system \eqref{eq:radical-abg} involves radical expressions, thus can not be directly handled by this method. However, the following substitutions could be helpful.
\begin{equation}\label{eq:subs-uv}
	u=\sqrt{\frac{x+z}{c_2}},~v=\sqrt{\frac{y+z}{c_1}}.
\end{equation}
We use $\bar a_0$, $\bar a_1$  and $\bar a_2$ to denote the results of $a_0$, $a_1$ and $a_2$ after the above substitutions, respectively. System \eqref{eq:radical-abg} is equivalent to the following polynomial system.
\begin{equation}\label{eq:poly-abg}
	\left\{\begin{split}
	&c_2u^2-(x+z)=0,\\
	&c_1v^2-(y+z)=0,\\
	&v-(x+y+z)=0,\\
	&u-(x+y+z)=0,\\
	&\numer\left(kz\left[-c_3+\frac{x+y}{(x+y+z)^2}\right]\right)=0,\\
	&\numer(1+\bar a_2+\bar a_1+\bar a_0)\cdot\denom(1+\bar a_2+\bar a_1+\bar a_0)>0,\\
	&\numer(1-\bar a_2+\bar a_1-\bar a_0)\cdot\denom(1-\bar a_2+\bar a_1-\bar a_0)>0,\\
	&\numer(-\bar a_0^2-\bar a_0\bar a_2+\bar a_1+1)\cdot\denom(-\bar a_0^2-\bar a_0\bar a_2+\bar a_1+1) >0,\\
	&\numer(-\bar a_0^2+\bar a_0\bar a_2-\bar a_1+1)\cdot\denom(-\bar a_0^2+\bar a_0\bar a_2-\bar a_1+1)>0,\\
	&x>0,~y>0,~z>0,~u>0,~v>0,\\
	&k>0,~1\geq l\geq 1,~c_1>0,~c_2>0,~c_3>0,
	\end{split}
	\right.
\end{equation}
where $\numer(\cdot)$ and $\denom(\cdot)$ stand for the numerator and the denominator.

In \cite{Tramontana2012H}, numerical simulations rather than symbolic deductions are used to explore the size of the parameter region for the stability of the $T_{ANB}$ model. The reason may be that the expressions of the eigenvalues of the Jacobian matrix are so complicated that any rigorous analysis seems to be impossible. However, by virtue of the Schur-Cohn criterion and the symbolic method in \cite{Li2014C}, we have obtained several analytical results.

It would be quite difficult to deal with the results of computations, if there exist $5$ parameters. Thus, as in \cite{Tramontana2012H}, we keep $c_3$ and $k$ as parameters and set $c_1=1.63$, $c_2=2.1$ and $l=0.6$ for the $T_{ANB}$ model, then we acquire the following theorem.

\begin{theorem}\label{thm:anb_stab}
Assume that $c_1=1.63$, $c_2=2.1$ and $l=0.6$. The equilibrium $E_2$ of the $T_{ANB}$ model is locally stable if and only if 
$$LS_{ANB}^1>0,~LS_{ANB}^2< 0,~LS_{ANB}^3>0~\text{and}~LS_{ANB}^4< 0,$$
where
\begin{equation}
LS_{ANB}^1 =100\,c_3-47,~LS_{ANB}^2=100\,c_3-373,
\end{equation}
and $LS_{ANB}^3$, $LS_{ANB}^4$ are complicated polynomials. See Appendix for details.
\end{theorem}

\begin{figure}[htbp]
    \centering
    \includegraphics[width=8cm]{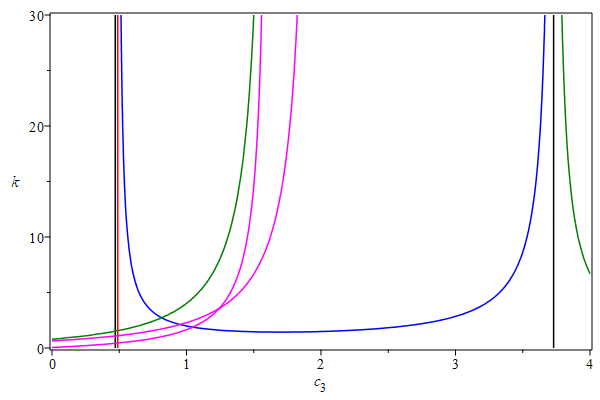}
    \caption{The $2$-dimensional $(c_3,k)$ parameter plane of the $T_{ANB}$ model with the other parameters fixed: $c_1=1.63$, $c_2=2.1$ and $l=0.6$. The two vertical black lines stand for $LS_{ANB}^1=0$ and $LS_{ANB}^2=0$. The blue curve is $LS_{ANB}^3=0$, while the green curve is  $LS_{ANB}^4=0$.}
    \label{fig:abg_1}
\end{figure}

The $(c_3,k)$ parameter plane of the $T_{ANB}$ model is depicted in Fig. \ref{fig:abg_1}. The equilibrium $E_2$ is locally stable in the region surrounded by the green curve, the blue curve, the two black vertical lines and the horizontal axis.

For the $T_{LNB}$ model, by setting $c_1=0.5$ and $c_2=0.55$ as in \cite{Naimzada2015T}, we obtain the following theorem.

\begin{theorem}\label{thm:lnb_stab}
Assume that $c_1=0.5$ and $c_2=0.55$. The equilibrium $E_2$ of the $T_{LNB}$ model is locally stable if and only if 
$$LS_{LNB}^1>0,~LS_{LNB}^2 < 0,~LS_{LNB}^3 < 0~\text{and}~LS_{LNB}^4 < 0,$$
where
\begin{equation}
LS_{LNB}^1 =20\,c_3-1,~LS_{LNB}^2 =20\,c_3-21,
\end{equation}
and the expressions of $LS_{LNB}^3$, $LS_{LNB}^4$ can be found in Appendix.
\end{theorem}

\begin{figure}[htbp]
    \centering
    \includegraphics[width=8cm]{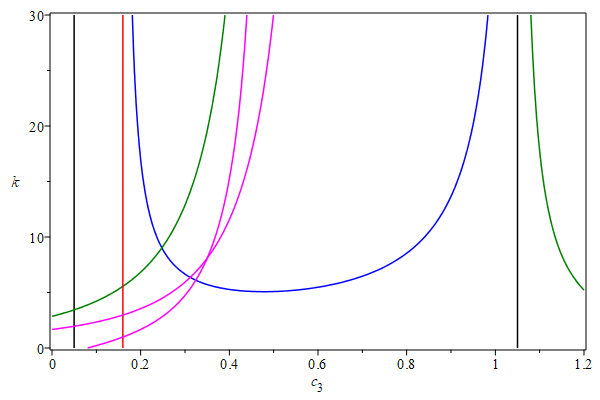}
    \caption{The $2$-dimensional $(c_3,k)$ parameter plane of the $T_{LNB}$ model with the other parameters fixed: $c_1=0.5$ and $c_2=0.55$. The two vertical black lines stand for $LS_{LNB}^1=0$ and $LS_{LNB}^2=0$. The blue curve is $LS_{LNB}^3=0$, while the green curve is $LS_{LNB}^4=0$.}
    \label{fig:lbg_1}
\end{figure}

The $(c_3,k)$ parameter plane of the $T_{LNB}$ model is given in Fig. \ref{fig:lbg_1}. The equilibrium $E_2$ is locally stable in the region surrounded by the green curve, the blue curve, the two black vertical lines and the horizontal axis.

Theorem \ref{thm:anb_stab} and Theorem \ref{thm:lnb_stab} confirm the observations  in \cite{Tramontana2012H,Naimzada2015T}. The reader may find that the stable regions of the two models have similar shapes. In Section  \ref{sec:pd-bifur} and \ref{sec:ns-bifur}, it would also be shown that the two models lose their stability through similar patterns.

Bischi and others \cite{Bischi2007O} proved that for a market with the isoelastic demand and an oligopoly of identical firms adopting the LMA mechanism, instability arises if there exist at least $5$ firms. While for an oligopoly of identical firms adopting the best response mechanism, the equilibrium goes unstable if at least $4$ firms compete with each other. This means that the LMA mechanism improves the stability compared to the best response mechanism for games with homogeneous players. In \cite{Naimzada2009C}, it was proved that the equilibrium is globally stable in a homogeneous duopoly game if both firms adopt the LMA mechanism. Moreover, Cavalli and Naimzada \cite{Cavalli2015N} investigated several heterogeneous duopoly games and found the similar behavior that the LMA mechanism may enlarge the size of the stable region, which suggests that stability is influenced not only by the rationality level of the firms but also by the reaction way the least rational firm uses to adapt its strategy. 

Therefore, it would be interesting to investigate how the LMA mechanism is connected to the stability for the heterogenous triopoly with the isoelastic demand. The comparison of the two games considered in this paper may be good material for this purpose. In order to make the results possible for comparison, we simplify the situation by setting the marginal costs of the three firms identical and acquire the following theorem.

\begin{theorem}\label{thm:const_c}
Suppose that $c_1=c_2=c_3=c$. For the $T_{ANB}$ model, the equilibrium $E_2$ is locally stable if and only if
\begin{equation}\label{eq:ckl}
	ck<\frac{192-102\,l}{68-33\,l}.
\end{equation}
Moreover, for the $T_{LNB}$ model, the equilibrium $E_2$ is locally stable if and only if
\begin{equation}\label{eq:ck}
	ck<\frac{62}{23}.
\end{equation}
\end{theorem}

\begin{figure}[htbp]
    \centering
    \includegraphics[width=8cm]{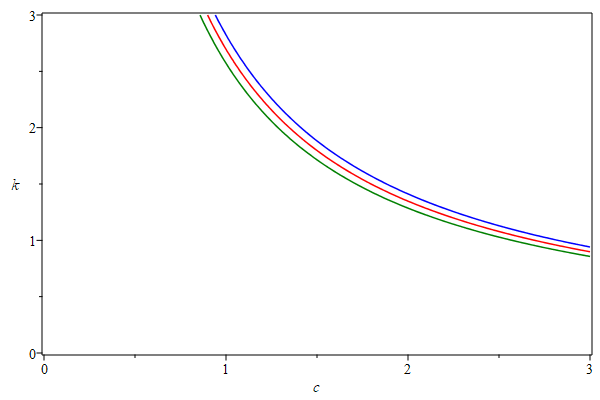}
    \caption{The $2$-dimensional $(c,k)$ parameter plane with identical marginal costs: $c_1=c_2=c_3=c$. }
    \label{fig:const_c}
\end{figure}

The parameter planes $(c,k)$ of the two games are depicted in Fig. \ref{fig:const_c}. The stable regions are surrounded by the horizontal axis, the vertical axis and the corresponding curves. The red curve represents the border of the region described by \eqref{eq:ck} as well as that by \eqref{eq:ckl} with $l=2/3$. The green and blue curves are the borders of \eqref{eq:ckl} with $l=1$ and $l=0$, respectively.

By Theorem \ref{thm:const_c}, it is easy to see that if $0\leq l<2/3$, the stable region of the $T_{ANB}$ model would be larger than that of the $T_{LNB}$ model. If $l=2/3$, the stable regions of these two models would be the same. Furthermore, if $2/3<l\leq 1$, the stable region of the $T_{ANB}$ model would be smaller than that of the $T_{LNB}$ model. As a special case, when $l=1$, the first firm in the $T_{ANB}$ model degenerates to a naive player adopting the best response mechanism. This confirms the observations in \cite{Bischi2007O,Naimzada2009C,Cavalli2015N} that the LMA mechanism has an effect of stability enhancement and suggests that we may extend a similar conclusion to the case of generic size heterogeneous oligopoly models with the isoelastic demand. However, if $l<1$, a general adaptive decision mechanism is used by the first firm in the $T_{ANB}$ model. As the proportion parameter $l$ goes from $1$ to $0$, the first firm adjusts its output with more caution and the stability region gets larger. It is also pointed out in \cite{Cavalli2015N} that the LMA mechanism might lead to a global underestimation of the prices and result in more careful adjustments. Thus, we conclude that the stability of an oligopoly may be connected essentially to the adjustment speed of the output.

\section{Period-doubling Bifurcation}\label{sec:pd-bifur}

The technique of linearization used for stability analysis may fail at bifurcation points, for near such points the dynamical behavior of the system may differ qualitatively from that of its linearization. This means that a small smooth change of the bifurcation parameters may cause a sudden qualitative change to the dynamic behavior of the system.

Similar to the continuous case, there are many different bifurcating situations for discrete systems. In this paper, we investigate two types of bifurcations by reducing the problem of bifurcation analysis to an algebraic problem.  In what follows, we first recall the formal definition of the period-doubling bifurcation, which could also be found in \cite{Kuznetsov1998E}.

\begin{definition}[Period-doubling Bifurcation]
	Consider a $n$-dimensional iteration map $F_{\uvar}:\,\xvar(t)\mapsto \xvar(t+1)$ with an equilibrium $E$, where $\xvar(t)=(x_1(t),x_2(t),\ldots,x_n(t))$ is the state vector at period $t$ and $\uvar=(u_1,u_2,\ldots,u_d)$ stands for all the parameters. A \emph{period-doubling bifurcation}, also called \emph{flip bifurcation}, takes place at a bifurcation parameter point $\uvar=\overline\uvar$ if and only if 
	the following two conditions are satisfied.
	\begin{enumerate}
		\item Eigenvalue assignment: the Jacobian matrix at $E$ has one real eigenvalue $\lambda_1(\overline\uvar)=1$, and $|\lambda_j(\overline\uvar)|<1$ for all the rest eigenvalues $\lambda_j$, $j\neq 1$.
		\item Transversality condition: for any parameter $u_j~(1\leq j \leq d)$,
		$$\frac{\partial |\lambda_1(\uvar)|}{\partial u_i}\Big|_{\uvar=\overline\uvar}\neq 0.$$
	\end{enumerate}
\end{definition}

From the economic point of view, there exist ranges of values of the parameters in which the dynamics are periodic if the period-doubling bifurcation occurs. It is quite realistic to assume that boundedly rational firms can not learn the pattern behind quantities and profits if long periods of periodic dynamics take place. 
In \cite{Wen2008A}, Wen and others generalized the Schur-Cohn criterion and proposed a new criterion for the detection of the period-doubling bifurcation. We revisit this criterion below.

\begin{proposition}
Let $A= \lambda^n + a_{n-1}\lambda^{n - 1} + \cdots + a_0$ be the characteristic polynomial of the Jacobian matrix of a $n$-dimensional iteration map $F_{\uvar}$ and $D^\pm_i$ be the same as in the Schur-Cohn criterion. A period-doubling bifurcation appears if and only if the following conditions are satisfied.
\begin{enumerate}
	\item Eigenvalue assignment: at the bifurcation point $\uvar=\overline\uvar$,
\begin{enumerate}
\item  $A(1)>0$ and $A(-1)=0$,
\item $D^\pm_1>0, D^\pm_3>0, \ldots, D^\pm_{n-3}>0, D^\pm_{n-1}>0$ (when $n$ is even), or\\[2pt] $D^\pm_2>0, D^\pm_4>0, \ldots, D^\pm_{n-3}>0, D^\pm_{n-1}>0$ (when $n$ is odd).
\end{enumerate}

	\item Transversality condition: for any parameter $u_j~(1\leq j\leq d)$,
	$$\frac{\sum_{i=0}^{n-1} (-1)^i \partial a_i/ \partial u_j}{\sum_{i=1}^{n-1}i(-1)^{n-i}a_i+n}\neq 0,$$
	where $\partial a_i/ \partial u_j$ is the partial derivative of $a_i(\uvar)$ with respective to $u_j$.
\end{enumerate}	
\end{proposition}

For the $3$-dimensional iteration maps \eqref{eq:iter-abg} and \eqref{eq:iter-lbg}, the criterion seems simple and is restated in the following corollary.

\begin{corollary}
	Consider a $3$-dimensional discrete dynamic system with $A=\lambda^3 + a_2\lambda^2 + a_1\lambda + a_0$ as the characteristic polynomial of its Jacobian matrix. There is a period-doubling bifurcation if and only if
	\begin{equation}
		\left\{\begin{split}
			&1+a_2+a_1+a_0>0,~-1+a_2-a_1+a_0=0,\\
			&-a_0^2-a_0a_2+a_1+1>0,~ -a_0^2+a_0a_2-a_1+1>0,\\
			&\frac{\partial a_0/\partial u_j-\partial a_1/ \partial u_j+\partial a_2/ \partial u_j}{3-2\,a_2+a_1}\neq 0,~\text{for any parameter}~u_j.
		\end{split}\right.
	\end{equation}
\end{corollary}

For the $T_{ANB}$ model, in order to calculate the bifurcation parameter surface determined by $-1+a_2-a_1+a_0=0$, we consider
\begin{equation}\label{eq:anb_pd}
	\left\{\begin{split}
	&-1+a_2-a_1+a_0=0,\\
	&\sqrt{\frac{y+z}{c_1}}-(x+y+z)=0,\\
	&\sqrt{\frac{x+z}{c_2}}-(x+y+z)=0,\\
	&kz\left[-c_3+\frac{x+y}{(x+y+z)^2}\right]=0,
	\end{split}
	\right.
\end{equation}
and then transform it to a polynomial system by applying the substitutions \eqref{eq:subs-uv}. Using the triangular decomposition method \cite{Wu1986B, Li2010D, Kalkbrener1993A, Wang2000C}, we obtain three triangular polynomial systems as follows. The union of the zeros of the resulting triangular polynomial systems is the same as the zero set of the system generated by substituting \eqref{eq:subs-uv} into \eqref{eq:anb_pd}.
\begin{equation}
\begin{split}
	\pset{T}_1=[~&u,~ v,~ z,~ y, ~x],\\
	\pset{T}_2=[~&u-z-y-x, ~v-z-y-x,~ c_3z+(-c_1+c_3)y-xc_1,\\
			&~ (-c_3+c_1-c_2)y+x(c_1-c_2+c_3),~ PD_{ANB}],\\
	\pset{T}_3=[~&(c_1+c_2)u-1, ~(c_1+c_2)v-1,~ z, ~(c_1^2+2\,c_1c_2+c_2^2)y-c_1, \\
	&(c_1^2+2\,c_1c_2+c_2^2)x-c_2,~ kc_3-c_1k-c_2k-2~],
\end{split}
\end{equation}
where $PD_{ANB}$ is complicated and can be found in Appendix. The polynomial $z$ is in $\pset{T}_1$ and $\pset{T}_3$, but the equilibrium $E_1$ is not of our concern. Thus, we focus only on $\pset{T}_2$ . It is clear that in $\pset{T}_2$ the last polynomial $PD_{ANB}$ involves only $c_1$, $c_2$, $c_3$, $l$ and $k$ (all the parameters) but none of $u$, $v$, $x$, $y$ and $z$, thus the bifurcation parameter surface of the period-doubling bifurcation for the $T_{ANB}$ model is $PD_{ANB}=0$. In addition, if the bifurcation could take place, it is also required the transversality condition is satisfied. Therefore, we obtain the following theorem.

\begin{theorem}
	For the $T_{ANB}$ model, the bifurcation parameter surface of the period-doubling bifurcation is $PD_{ANB}=0$. Furthermore, if a period-doubling bifurcation takes place, the following transversality condition should be satisfied.
	$$PT_{ANB}^k\neq 0,~ PT_{ANB}^l\neq 0, ~PT_{ANB}^{c_1}\neq 0, ~PT_{ANB}^{c_2}\neq 0, ~PT_{ANB}^{c_3}\neq 0.$$
	See Appendix for details.
\end{theorem}

\begin{figure}[htbp]
    \centering
    \includegraphics[width=8cm]{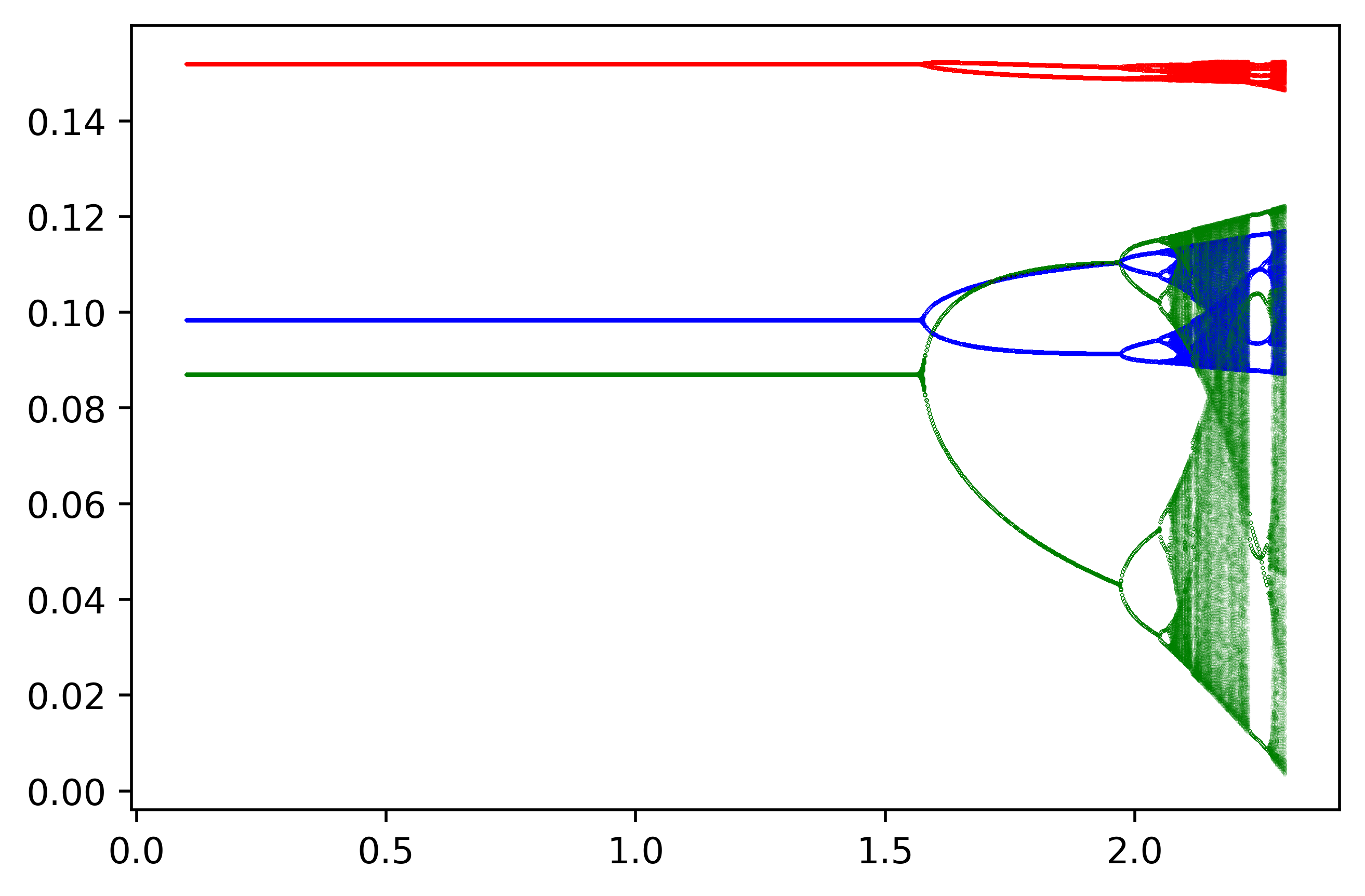}
    \caption{The bifurcation diagram with respect to $k$ for the $T_{ANB}$ model with $c_1=1.63$, $c_2=2.1$, $c_3=2.2$ and $l=0.6$. }
    \label{fig:bif_abg_pd}
\end{figure}

Similarly, we have the following theorem for the $T_{LNB}$ model.

\begin{theorem}
	For the $T_{LNB}$ model, the bifurcation parameter surface of the period-doubling bifurcation is $PD_{LNB}=0$. Furthermore, if a period-doubling bifurcation takes place, the following transversality condition should be satisfied.
	$$PT_{LNB}^k\neq 0, ~PT_{LNB}^{c_1}\neq 0, ~PT_{LNB}^{c_2}\neq 0,~ PT_{LNB}^{c_3}\neq 0.$$
	See Appendix for details.
\end{theorem}

\begin{figure}[htbp]
    \centering
    \includegraphics[width=8cm]{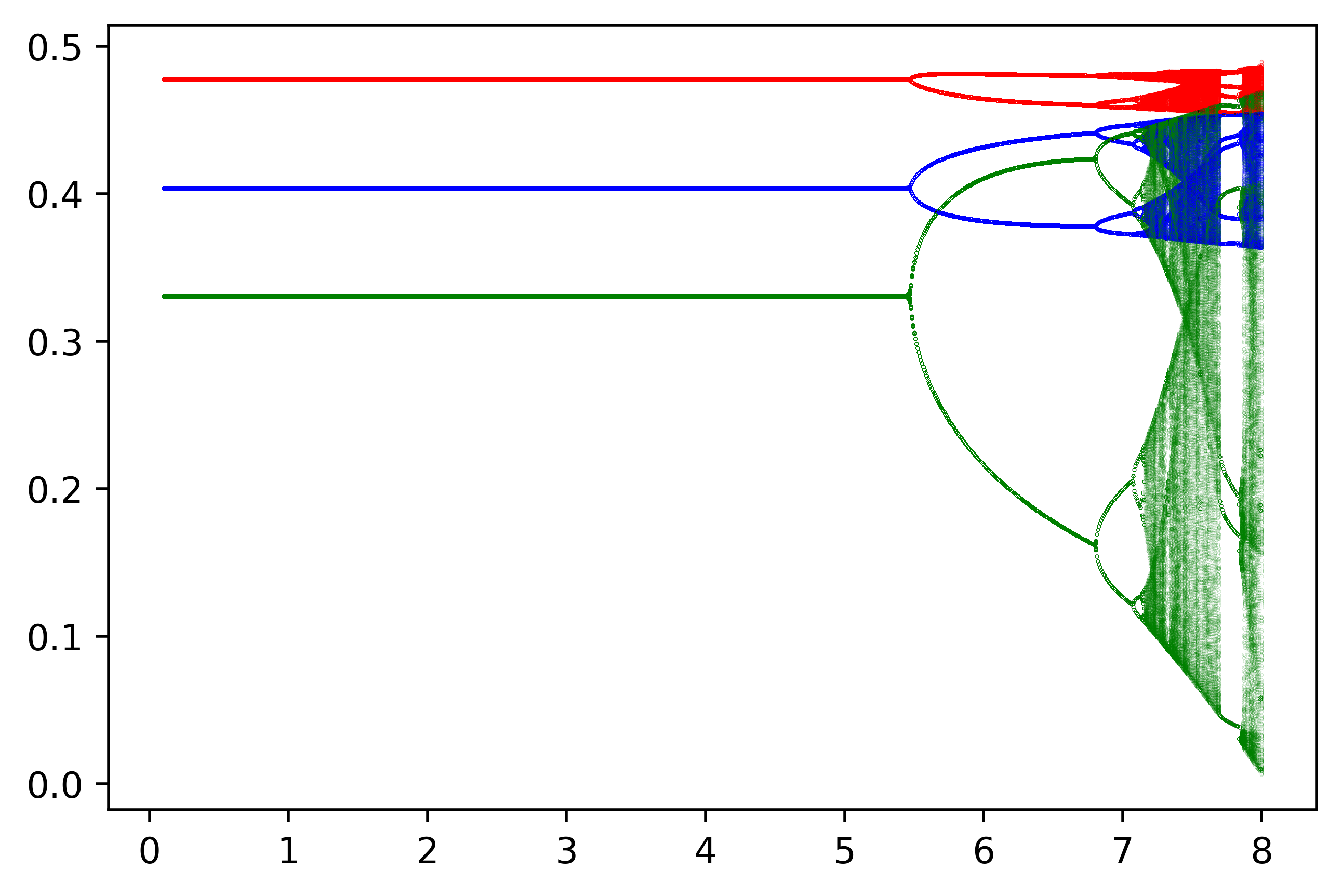}
    \caption{The bifurcation diagram with respect to $k$ for the $T_{LNB}$ model with $c_1=0.5$, $c_2=0.55$ and $c_3=0.6$. }
    \label{fig:bif_lbg_pd}
\end{figure}

\begin{figure}[htbp]
    \centering
    \includegraphics[width=8cm]{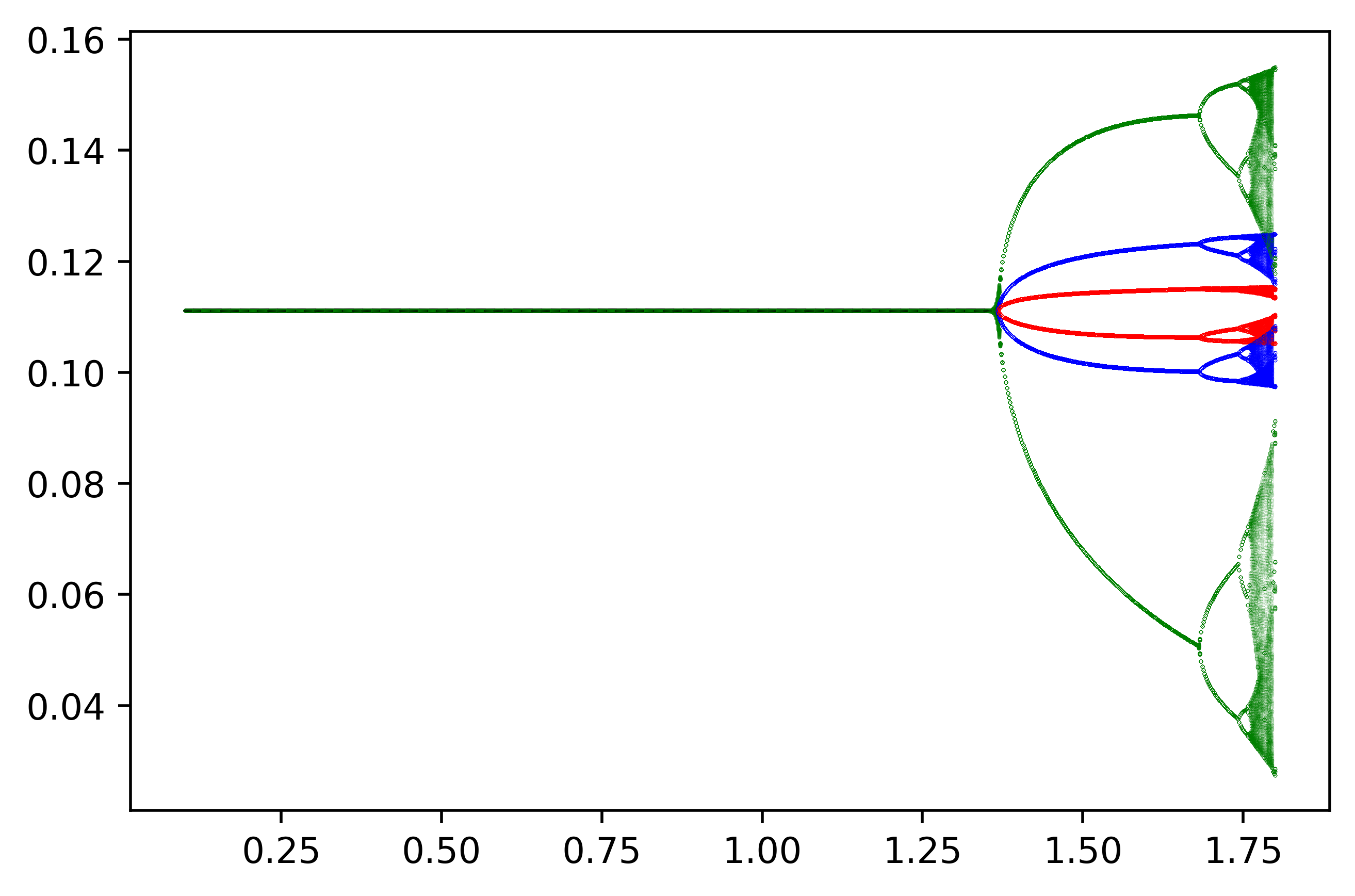}
    \caption{The bifurcation diagram with respect to $k$ for the $T_{ANB}$ model with identical marginal costs $c_1=c_2=c_3=2$ and $l=0.5$. }
    \label{fig:bif_abg_c}
\end{figure}

\begin{figure}[htbp]
    \centering
    \includegraphics[width=8cm]{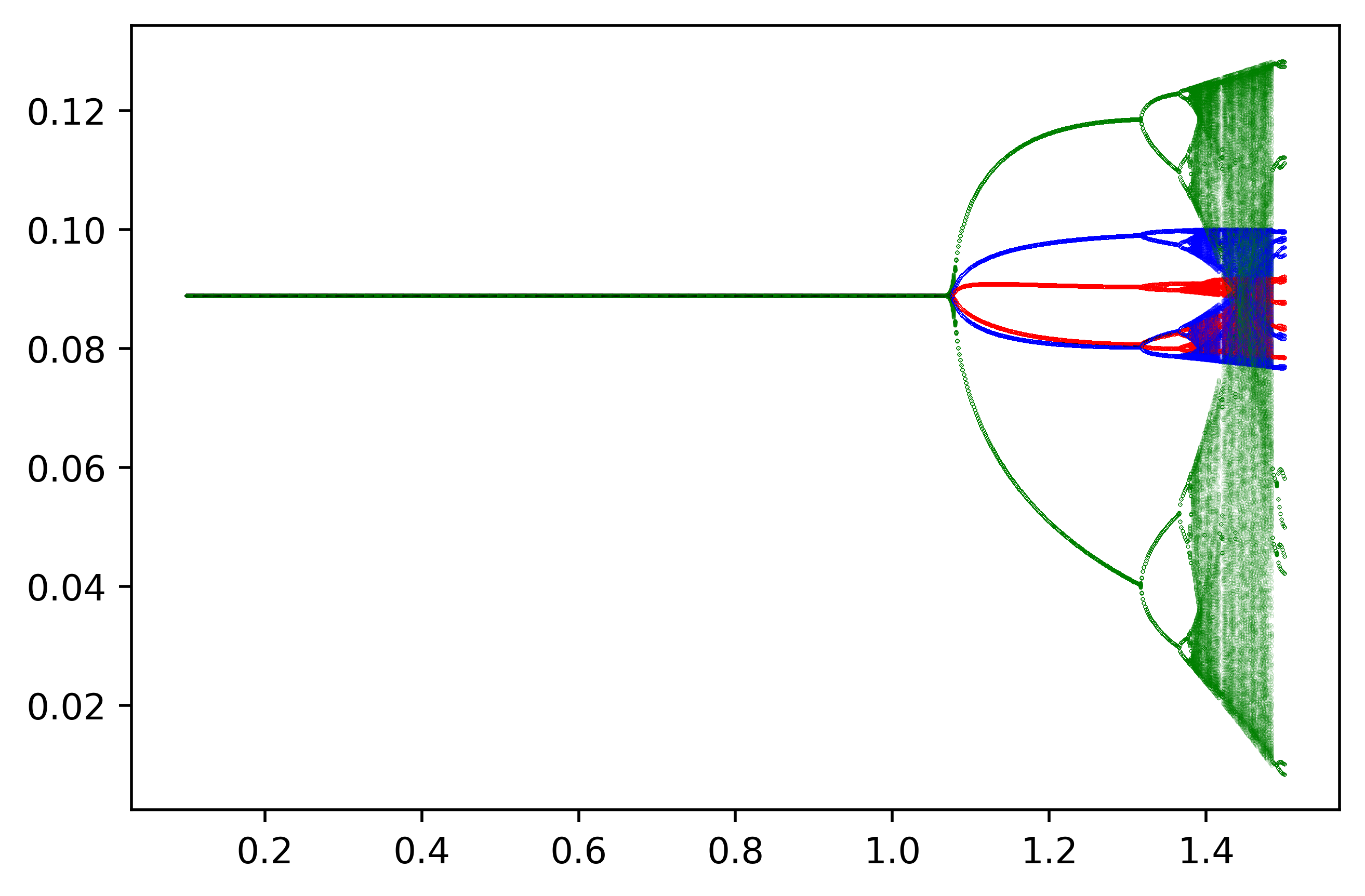}
    \caption{The bifurcation diagram with respect to $k$ for the $T_{LNB}$ model with identical marginal costs $c_1=c_2=c_3=2.5$. }
    \label{fig:bif_lbg_c}
\end{figure}

It is easy to verify that $PD_{ANB}$ becomes $LS_{ANB}^3$ if setting $c_1=1.63$, $c_2=2.1$, $l=0.6$ as in Theorem \ref{thm:anb_stab}, and $PD_{LNB}$  becomes $LS_{LNB}^3$ if setting $c_1=0.5$, $c_2=0.55$ as in Theorem \ref{thm:lnb_stab}. This proves that the equilibrium $E_2$ would lose its stability through the period-doubling bifurcation when a parameter point in the stable region goes across the blue curve of Fig. \ref{fig:abg_1} or Fig. \ref{fig:lbg_1}. Furthermore, for we keep only $c_3$ and $k$ as the parameters in Section \ref{sec:stability}, the transversality conditions with respect to $c_3$ and $k$ are required, which are $PT^k_{ANB}\neq 0, PT^{c_3}_{ANB}\neq 0$ for the $T_{ANB}$ model and  $PT^k_{LNB}\neq 0, PT^{c_3}_{LNB}\neq 0$ for the $T_{LNB}$ model, respectively. These conditions are shown in Fig. \ref{fig:abg_1} and Fig. \ref{fig:lbg_1} with the red curves. It is worth noting that the red curves do not pass through the segments of the blue curves where the period-doubling bifurcation takes place for both models. Therefore, the problem that the transversality condition is violated does not occur for our settings.

Furthermore, for the case of the identical marginal costs, $PD_{ANB}$ would be $(33\,ckl-68\,ck-102\,l+192)c^3$ corresponding to \eqref{eq:ckl}, and $PD_{ANB}$ would be $(23\,ck-62)c^2$ corresponding to \eqref{eq:ck}. This means that for both models the equilibrium $E_2$ loses its stability via the period-doubling bifurcation if assuming that the marginal costs of the firms are identical. The bifurcation diagrams with respect to $k$ are given in Fig. \ref{fig:bif_abg_c} and  Fig. \ref{fig:bif_lbg_c}.

\section{Neimark-Sacker Bifurcation}\label{sec:ns-bifur}

The {Neimark-Sacker bifurcation} for discrete dynamical systems corresponds to the Hopf bifurcation in the continuous case. When a Neimark–Sacker bifurcation takes place, dynamics suddenly become quasiperiodic, which is more difficult to deal with than a periodic solution. This is because that a quasiperiodic trajectory looks like chaos to boundedly rational firms with incomplete information and limited capability of calculation. The reader may refer to \cite{Guckenheimer2013N} for the following formal definition of the Neimark-Sacker bifurcation.

\begin{definition}[Neimark-Sacker Bifurcation]
Consider a $n$-dimensional iteration map $F_{\uvar}:\,\xvar(t)\mapsto \xvar(t+1)$ with an equilibrium $E$, where $\xvar(t)=(x_1(t),x_2(t),\ldots,x_n(t))$ is the state vector at period $t$ and $\uvar=(u_1,u_2,\ldots,u_d)$ stands for all the parameters. A \emph{Neimark-Sacker bifurcation} takes place at a bifurcation parameter point $\uvar=\overline\uvar$ if and only if 
	the following two conditions are satisfied.
	\begin{enumerate}
		\item Eigenvalue assignment: the Jacobian matrix at $E$ has a pair of complex conjugate eigenvalues $\lambda_1(\overline\uvar)$, $\lambda_2(\overline\uvar)$ on the unit circle, i.e., $|\lambda_1(\overline\uvar)|=|\lambda_2(\overline\uvar)|=1$, and $|\lambda_j(\overline\uvar)|<1$ for all the rest eigenvalues $\lambda_j$, $j\neq 1,2$.
		\item Transversality condition: for any parameter $u_j~(1\leq j \leq d)$,
		$$\frac{\partial |\lambda_1(\uvar)|}{\partial u_i}\Big|_{\uvar=\overline\uvar}\neq 0.$$
	\end{enumerate}

\end{definition}

Based on the Schur-Cohn criterion, a criterion for the identification of the Neimark-Sacker bifurcation is proposed in \cite{Wen2005C}, which is restated below.

\begin{proposition}
Let $A= \lambda^n + a_{n-1}\lambda^{n - 1} + \cdots + a_0$ be the characteristic polynomial of the Jacobian matrix of a $n$-dimensional iteration map $F_{\uvar}$ and $D^\pm_i$ be the same as in the Schur-Cohn criterion. A Neimark-Sacker bifurcation occurs if and only if the following conditions are satisfied.
	\begin{enumerate}
		\item Eigenvalue assignment: at the bifurcation point $\uvar=\overline\uvar$,
			\begin{enumerate}
			\item $A(1)>0$ and $(-1)^nA(-1)>0$,
			
			\item $D^\pm_1>0, D^\pm_3>0, \ldots, D^\pm_{n-3}>0, D^+_{n-1}>0, D^-_{n-1}=0$
			(when $n$ is even), or\\[2pt]
			$D^\pm_2>0, D^\pm_4>0, \ldots, D^\pm_{n-3}>0, D^+_{n-1}>0,
			D^-_{n-1}=0$ (when $n$ is odd).
			\end{enumerate}
		
		\item Transversality condition: for any parameter $u_j~(1\leq j\leq d)$,
		$$\frac{\partial D^-_{n-1}}{\partial u_j}\Big|_{\uvar=\overline\uvar}\neq 0.$$
	\end{enumerate}

\end{proposition}

\begin{corollary}
	For a $3$-dimensional discrete dynamic system with $A=\lambda^3 + a_2\lambda^2 + a_1\lambda + a_0$ as the characteristic polynomial of its Jacobian matrix, the Neimark-Sacker bifurcation occurs if and only if
	\begin{equation}
		\left\{\begin{split}
			&1+a_2+a_1+a_0>0,~-1+a_2-a_1+a_0<0,\\
			&-a_0^2-a_0a_2+a_1+1>0,~ -a_0^2+a_0a_2-a_1+1=0,\\
			&-2\frac{\partial a_0}{\partial u_j}+a_0\frac{\partial a_2}{\partial u_j}+a_2\frac{\partial a_0}{\partial u_j}-\frac{\partial a_1}{\partial u_j}\neq 0,~\text{for any parameter}~u_j.
		\end{split}\right.
	\end{equation}
\end{corollary}

Based on a series of computations,  the following two theorems are obtained.

\begin{theorem}
	For the $T_{ANB}$ model, the bifurcation parameter surface of the Neimark-Sacker bifurcation is $NS_{ANB}=0$. Furthermore, if a Neimark-Sacker bifurcation appears, the following transversality condition should be satisfied.
	$$NT_{ANB}^k\neq 0,~ NT_{ANB}^l\neq 0,~ NT_{ANB}^{c_1}\neq 0,~ NT_{ANB}^{c_2}\neq 0,~ NT_{ANB}^{c_3}\neq 0.$$
	See Appendix for details.
\end{theorem}

\begin{theorem}
	For the $T_{LNB}$ model, the bifurcation parameter surface of the Neimark-Sacker bifurcation is $NS_{LNB}=0$. Furthermore, if a Neimark-Sacker bifurcation occurs, the following transversality condition should be satisfied.
	$$NT_{LNB}^k\neq 0,~ NT_{LNB}^{c_1}\neq 0,~ NT_{LNB}^{c_2}\neq 0,~ NT_{LNB}^{c_3}\neq 0.$$
	See Appendix for details.
\end{theorem}

\begin{figure}[htbp]
    \centering
    \includegraphics[width=8cm]{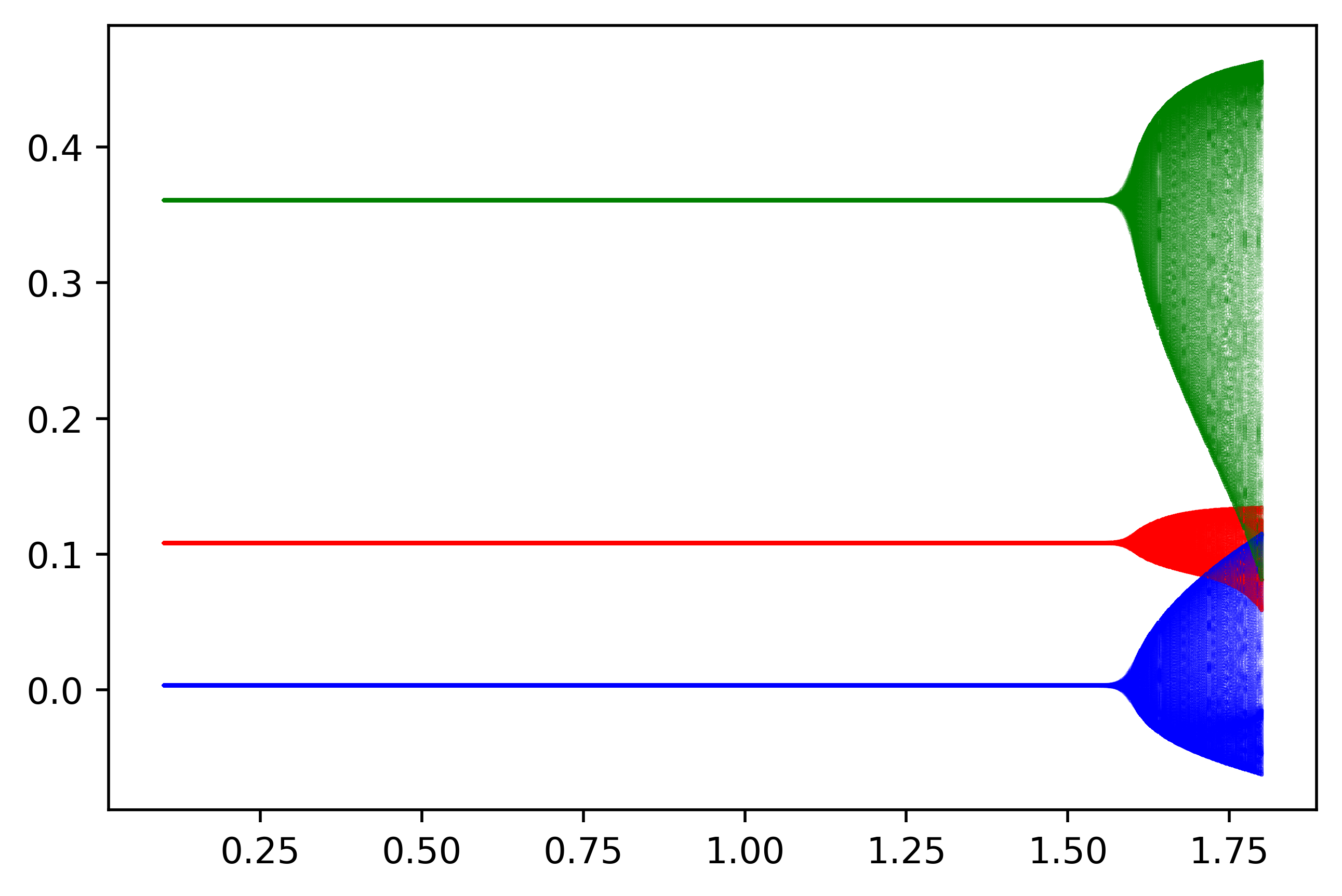}
    \caption{The bifurcation diagram with respect to $k$ for the $T_{ANB}$ model with $c_1=1.63$, $c_2=2.1$, $c_3=0.5$ and $l=0.6$. }
    \label{fig:bif_abg_ns}
\end{figure}

\begin{figure}[htbp]
    \centering
    \includegraphics[width=8cm]{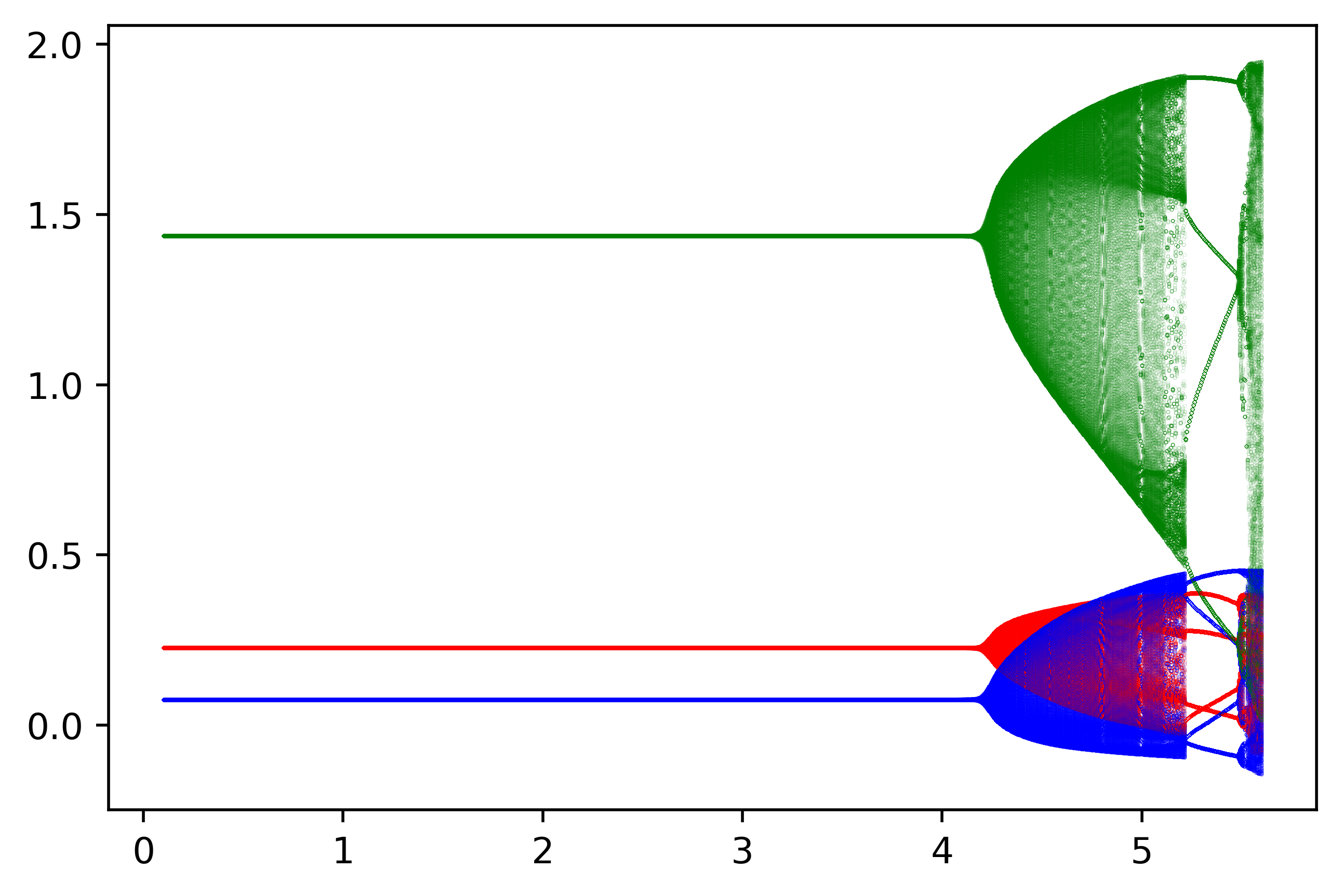}
    \caption{The bifurcation diagram with respect to $k$ for the $T_{LNB}$ model with $c_1=0.5$, $c_2=0.55$ and $c_3=0.1$. }
    \label{fig:bif_lbg_ns}
\end{figure}

We could see that $NS_{ANB}$ and $NS_{LNB}$ would be $LS_{ANB}^4$ and $LS_{LNB}^4$, respectively, if setting $c_1=1.63$, $c_2=2.1$, $l=0.6$ as in Theorem \ref{thm:anb_stab} and setting $c_1=0.5$, $c_2=0.55$ as in Theorem \ref{thm:lnb_stab}. This means that across the green curve in Fig. \ref{fig:abg_1} or Fig. \ref{fig:lbg_1} as the parameters vary, the equilibrium $E_2$ would lose its stability through the Neimark-Sacker bifurcation. Moreover, as we keep only $c_3$ and $k$ as the parameters in these two figures, the transversality conditions with respect to $c_3$ and $k$ are needed, which are $NT^k_{ANB}\neq 0, NT^{c_3}_{ANB}\neq 0$ for the $T_{ANB}$ model and  $NT^k_{LNB}\neq 0, NT^{c_3}_{LNB}\neq 0$ for the $T_{LNB}$ model. These conditions are shown in Fig. \ref{fig:abg_1} and Fig. \ref{fig:lbg_1} with magenta curves. It is obvious that the magenta curves do not pass through the segments of the green curves where the Neimark-Sacker bifurcation is developed for both models. Therefore, for our settings, we do not need to worry about the possibility that the transversality condition might be violated.

\section{Concluding Remarks}\label{sec:con}

In this paper, we have analyzed dynamic behaviors of two different heterogeneous triopoly games under the assumption that the demand function of the market is isoelastic. Nonlinearities are introduced both in the demand function and in the decisional mechanism adopted by the firms. We have considered several decisional mechanisms including the naive (best response) mechanism, the boundedly rational (gradient adjustment) mechanism, the adaptive mechanism and the local monopolistic approximation (LMA) mechanism. In particular, in both games, the second firm is a naive player, while the third firm is a boundedly rational player. The only difference is that firm 1 is an adaptive player in the $T_{ANB}$ model, but an LMA player in the $T_{LNB}$ model.

The existing results of triopoly games in the literature are usually based on observations through numerical simulations. By contrast, using the approach proposed in \cite{Li2014C}, we have obtained several analytical results for the two models considered in this paper. In particular, we have rigorously proved the existence of double routes, through the period-doubling bifurcation and through the Neimark-Sacker bifurcation, from the stable state to the chaotic behavior. In addition, we have considered the special case that the involved firms have identical marginal costs and acquired the necessary and sufficient condition of the local stability for each model. It has been found that the presence of the local monopolistic approximation (LMA) mechanism might have a stabilizing effect for heterogeneous triopoly games with the isoelastic demand.

As our future work, the rigorous proof of the existence of chaos is theoretically interesting although it might have been observed via simulations. Moreover, we are also curious to learn how to control chaos in oligopoly games using symbolic methods.

\section*{Acknowledgments}

The author is grateful to the anonymous referees for their helpful comments. This work has been supported by Major Research and Cultivation Project of Dongguan City College (No. 2021YZDYB04Z) and Social Development Science and Technology Project of Dongguan (No. 20211800900692).

\bibliographystyle{abbrv}
\bibliography{ref_LMA1}

\section*{Appendix}

\begin{align*}

\begin{autobreak}
LS_{ANB}^3 =
300000000\,c_3^{4} k
-15102000000\,c_3^{3} k
-600000000\,c_3^{3}
+196817000000\,c_3^{2} k
+2238000000\,c_3^{2}
-607025944200\,c_3k
+85553100000\,c_3
+251560953577\,k
+256838922600,
\end{autobreak}\\

\begin{autobreak}
LS_{ANB}^4 =
900000000000000\,c_3^{7} k^{2}
-41949000000000000\,c_3^{6} k^{2}
-1800000000000000\,c_3^{6} k
+667715010000000000\,c_3^{5} k^{2}
+52020000000000000\,c_3^{5} k
-4657740430500000000\,c_3^{4} k^{2}
+ 900000000000000\,c_3^{5}
-523900860000000000\,c_3^{4} k
+14725627073035000000\,c_3^{3} k^{2}
-10071000000000000\,c_3^{4}
+1398076600800000000\,c_3^{3} k
-15544953327461830000\,c_3^{2} k^{2}
+56313540000000000\,c_3^{3} 
+7799709933666000000\,c_3^{2} k
-13291730300598585700\,c_3k^{2}
-23227083000000000\,c_3^{2}
-29520748208318520000\,c_3k
+25092830446277361973\,k^{2}
-5485330423251000000\,c_3
-2258144963547738600\,k
- 14615486984757870000,
\end{autobreak}\\

\begin{autobreak}
LS_{LNB}^3 =
40000\,c_3^{3} k
-218800\,c_3^{2} k
+16000\,c_3^{2}
+212860\,c_3k
-174400\,c_3
-28581\,k
-59960,	
\end{autobreak}\\

\begin{autobreak}
LS_{LNB}^4 =
3200000\,c_3^{6} k^{2}
-3200000\,c_3^{5} k^{2}
+3200000\,c_3^{5} k
-14704000\,c_3^{4} k^{2}
-4800000\,c_3^{4} k
+28668800\,c_3^{3} k^{2}
+800000\,c_3^{4}
+11776000\,c_3^{3} k
-15143500\,c_3^{2} k^{2}
+160000\,c_3^{3}
+11631200\,c_3^{2}k
+9660\,c_3k^{2}
-1940000\,c_3^{2}
-25391060\,c_3k
+1166445\,k^{2}
-15682800\,c_3
+1955352\,k
-15115275,
\end{autobreak}\\

\begin{autobreak}
PD_{ANB} = 
k l c_3^{4}
+(2\,c_1k l
-10\,c_2k l
-12\,c_1k
-2\,l) c_3^{3}
+(
-12\,c_1c_2k l
+4\,c_2^{2} k l
+4\,c_1^{2} k
+116\,c_1c_2k
+2\,c_1l
+2\,c_2l) c_3^{2} 
+(
-2\,c_1^{3} k l
+6\,c_1^{2} c_2k l
+18\,c_1c_2^{2} k l
+10\,c_2^{3} k l
+12\,c_1^{3} k
-104\,c_1^{2} c_2k
-116\,c_1c_2^{2} k
+10\,c_1^{2} l 
-44\,c_1c_2l
+10\,c_2^{2} l
+64\,c_1c_2) c_3
-c_1^{4} k l
+8\,c_1^{3} c_2k l
+14\,c_1^{2} c_2^{2} k l
-5\,c_2^{4} k l 
-4\,c_1^{4} k
+4\,c_1^{3} c_2k
+20\,c_1^{2} c_2^{2} k
+12\,c_1c_2^{3} k
+6\,c_1^{3} l
-46\,c_1^{2} c_2l
-46\,c_1c_2^{2} l
+6\,c_2^{3} l
+64\,c_1^{2} c_2
+64\,c_1c_2^{2},
\end{autobreak}\\

\begin{autobreak}
PT^k_{ANB}=
l c_3^{3}
+(3\,c_1l
-9\,c_2l
-12\,c_1) c_3^{2}
+(3\,c_1^{2} l
-18\,c_1c_2l
-5\,c_2^{2} l
-8\,c_1^{2}
+104\,c_1c_2) c_3
+c_1^{3} l
-9\,c_1^{2} c_2l
-5\,c_1c_2^{2} l
+5\,c_2^{3} l
+4\,c_1^{3}
-8\,c_1^{2} c_2
-12\,c_1c_2^{2},
\end{autobreak}\\

\begin{autobreak}
PT^l_{ANB}=
k c_3^{3}
+(c_1k
-11\,c_2k
-2) c_3^{2}
+(
-c_1^{2} k
-2\,c_1c_2k
+15\,c_2^{2} k
+4\,c_1
+4\,c_2) c_3
-c_1^{3} k
+9\,c_1^{2} c_2k
+5\,c_1c_2^{2} k
-5\,c_2^{3} k
+6\,c_1^{2}
-52\,c_1c_2
+6\,c_2^{2},
\end{autobreak}\\

\begin{autobreak}
PT^{c_1}_{ANB}=
k c_3^{3}
+(
-c_1k
-11\,c_2k
-2) c_3^{2}
+(
-c_1^{2} k
+14\,c_1c_2k
+15\,c_2^{2} k
-4\,c_1
+4\,c_2) c_3
+c_1^{3} k
-3\,c_1^{2} c_2k
-9\,c_1c_2^{2} k
-5\,c_2^{3} k
-2\,c_1^{2}
+4\,c_1c_2
+6\,c_2^{2},
\end{autobreak}\\

\begin{autobreak}
PT^{c_2}_{ANB}=
(k l
+4\,k) c_1^{3}
+(c_2k l
+c_3k l
+8\,c_2k
-16\,c_3k
-6\,l) c_1^{2}
+(
-c_2^{2} k l
+2\,c_2c_3k l
-c_3^{2} k l
+4\,c_2^{2} k
-16\,c_2c_3  k
+12\,c_3^{2} k
-4\,c_2l
-4\,c_3l) c_1
-c_2^{3} k l
+c_2^{2} c_3k l
+c_2\,c_3^{2} k l
-c_3^{3} k l
+2\,c_2^{2} l
+4\,c_2c_3l
+2\,c_3^{2} l,
\end{autobreak}\\

\begin{autobreak}
PT^{c_3}_{ANB}=
l c_3^{2}
+(
-2\,c_1l
-2\,c_2l) c_3
-3\,c_1^{2} l
+26\,c_1c_2l
-3\,c_2^{2} l
-32\,c_1c_2,
\end{autobreak}\\

\begin{autobreak}
PD_{LNB}=
(3\,c_1^{3}
-11\,c_1^{2} c_2
-5\,c_1^{2} c_3
-15\,c_1c_2^{2}
+54\,c_1c_2c_3
-3\,c_1c_3^{2}
-c_2^{3}
+43\,c_2^{2} c_3
-47\,c_2c_3^{2}
+5\,c_3^{3}) k
-6\,c_1^{2}
+20\,c_1c_2
-4\,c_1c_3
-38\,c_2^{2}
-36\,c_2c_3
+2\,c_3^{2},
\end{autobreak}\\

\begin{autobreak}
PT^{k}_{LNB}=
5\,c_3^{2}
+(2\,c_1
-42\,c_2) c_3
-3\,c_1^{2}
+14\,c_1c_2
+c_2^{2},
\end{autobreak}\\

\begin{autobreak}
PT^{c_1}_{LNB}=
(k c_3
-c_1k
-c_2k
+4) (c_3
+c_1
-7\,c_2) ,
\end{autobreak}\\

\begin{autobreak}
PT^{c_2}_{LNB}=
5\,k c_3^{2}
+(
-8\,c_1k
-6\,c_2k
+2) c_3
+3\,c_1^{2} k
+4\,c_1c_2k
+c_2^{2} k
-6\,c_1
+2\,c_2,
\end{autobreak}\\

\begin{autobreak}
PT^{c_3}_{LNB}=
c_3^{2}
+(
-2\,c_1
-18\,c_2) c_3
-3\,c_1^{2}
+10\,c_1c_2
-19\,c_2^{2},
\end{autobreak}\\

\begin{autobreak}
NS_{ANB} =
k^{2} l^{2} c_3^{7}
+(3\,c_1k^{2} l^{2}
-9\,c_2k^{2} l^{2}
-12\,c_1k^{2} l
-2\,k l^{2}) c_3^{6}
+(c_1^{2} k^{2} l^{2}
-22\,c_1c_2k^{2} l^{2}
+29\,c_2^{2} k^{2} l^{2}
-8\,c_1^{2} k^{2} l
+80\,c_1c_2k^{2} l
+36\,c_1^{2} k^{2}
+12\,c_2k l^{2}
+12\,c_1k l
+l^{2}) c_3^{5}
+(
-5\,c_1^{3} k^{2} l^{2}
-3\,c_1^{2} c_2k ^{2} l^{2}
+61\,c_1c_2^{2} k^{2} l^{2}
-37\,c_2^{3} k^{2} l^{2}
+28\,c_1^{3} k^{2} l
+16\,c_1^{2} c_2k^{2} l
-204\,c_1c_2^{2} k^{2} l
-60\,c_1^{3} k^{2}
-156\,c_1^{2} c_2k^{2}
+14\,c_1^{2} k l^{2}
-14\,c_2^{2} k l^{2}
-28\,c_1^{2} k l
-76\,c_1c_2k l
-3\,c_1l^{2}
-3\,c_2l^{2}) c_3^{4}
+(
-5\,c_1^{4} k^{2} l^{2}
+28\,c_1^{3} c_2k^{2} l^{2}
+10\,c_1^{2} c_2^{2} k^{2} l^{2}
-84\,c_1\,c_2^{3} k^{2} l^{2}
+3\,c_2^{4} k^{2} l^{2}
+16\,c_1^{4} k^{2} l
-128\,c_1^{3} c_2k^{2} l
+16\,c_1^{2} c_2^{2} k^{2} l
+288\,c_1c_2^{3} k^{2} l
-8\,  c_1^{4} k^{2}
+304\,c_1^{3} c_2k^{2}
+16\,c_1^{3} k l^{2}
+120\,c_1^{2} c_2^{2} k^{2}
-56\,c_1^{2} c_2k l^{2}
-24\,c_2^{3} k l^{2}
-40\,c_1^{3} k l
+16\,c_1^{2} c_2k l
+56\,c_1c_2^{2} k l
+96\,c_1^{2} c_2k
-6\,c_1^{2} l^{2}
+4\,c_1c_2l^{2}
-6\,c_2^{2} l^{2}
+16\,  c_1c_2l) c_3^{3}
+(c_1^{5} k^{2} l^{2}
+17\,c_1^{4} c_2k^{2} l^{2}
-46\,c_1^{3} c_2^{2} k^{2} l^{2}
-38\,c_1^{2} c_2^{3} k^{2} l^{2}
+61\,c_1c_2^{4} k^{2} l^{2}
+37\,c_2^{5} k^{2} l^{2}
-20\,c_1^{5} k^{2} l
+184\,c_1^{3} c_2^{2} k^{2} l
-96\,c_1^{2} c_2^{3} k^{2} l
- 260\,c_1c_2^{4} k^{2} l
+56\,c_1^{5} k^{2}
-216\,c_1^{4} c_2k^{2}
-6\,c_1^{4} k l^{2}
-152\,c_1^{3} c_2^{2} k^{2}
-48\,c_1^{3} c_2k l^{2}
+120\,c_1^{2} c_2^{3} k^{2}
+52\,c_1^{2} c_2^{2} k l^{2}
+34\,c_2^{4} k l^{2}
+40\,c_1^{4} k l
+56\,c_1^{3} c_2k l
+248\,c_1^{2} c_2^{2} k l
-24\,c_1c_2^{3} k l
-32\,c_1^{3} c_2k
+10\,c_1^{3} l^{2}
+96\,c_1^{2} c_2^{2} k
+14\,c_1^{2} c_2l^{2}
+14\,c_1c_2^{2} l^{2}
+10\,c_2^{3} l^{2}
-16\,c_1^{2} c_2l
-16\,c_1c_2^{2} l) c_3^{2}
+(3\,c_1^{6} k^{2} l^{2}
-6\,c_1^{5} c_2k^{2} l^{2}
-23\,c_1 ^{4} c_2^{2} k^{2} l^{2}
+28\,c_1^{3} c_2^{3} k^{2} l^{2}
+53\,c_1^{2} c_2^{4} k^{2} l^{2}
-22\,c_1c_2^{5} k^{2} l^{2}
-33\,c_2^{6} k^{2} l^{2}
-8\,c_1^{6} k^{2} l
+48\,c_1^{5} c_2k^{2} l
-48\,c_1^{4} c_2^{2} k^{2} l
-128\,c_1^{3} c_2^{3} k^{2} l
+120\,c_1^{2} c_2^{4} k^{2}  l
+144\,c_1c_2^{5} k^{2} l
-28\,c_1^{6} k^{2}
+80\,c_1^{5} c_2k^{2}
-16\,c_1^{5} k l^{2}
+88\,c_1^{4} c_2^{2} k^{2}
+12\,c_1^{4} c_2k l^{2}
-176\,c_1^{3} c_2^{3} k^{2}
+48\,c_1^{3} c_2^{2} k l^{2}
-156\,c_1^{2} c_2^{4} k^{2}
+8\,c_1^{2} c_2^{3} k l^{2}
+12\,c_2^{5} k  l^{2}
+28\,c_1^{5} k l
+16\,c_1^{4} c_2k l
-120\,c_1^{3} c_2^{2} k l
-176\,c_1^{2} c_2^{3} k l
-68\,c_1c_2^{4} k l
-96\,c_1^{4} c_2k
+21\,c_1^{4} l^{2}
-192\,c_1^{3} c_2^{2} k
+4\,c_1^{3} c_2l^{2}
-96\,c_1^{2} c_2^{3} k
-34\,c_1^{2} c_2^{2} l^{2}
+4\,c_1c_2^{3}  l^{2}
+21\,c_2^{4} l^{2}
-80\,c_1^{3} c_2l
-160\,c_1^{2} c_2^{2} l
-80\,c_1c_2^{3} l) c_3
+9\,c_1^{5} l^{2}
+52\,c_1^{5} c_2k l
+72\,c_1^{4} c_2^{2} k l
+40\,c_1^{3} c_2^{3} k l
+132\,c_1^{2} c_2^{4} k l
+100\,c_1c_2^{5} k l
-5\,c_1^{6} c_2k^{2} l^{2}
+c_1^{5}  c_2^{2} k^{2} l^{2}
+19\,c_1^{4} c_2^{3} k^{2} l^{2}
-5\,c_1^{3} c_2^{4} k^{2} l^{2}
-23\,c_1^{2} c_2^{5} k^{2} l^{2}
+3\,c_1c_2^{6} k^{2} l^{2}
-16\,c_1^{6} c_2k^{2} l
-12\,c_1^{5} c_2^{2} k^{2} l
+64\,c_1^{4} c_2^{3} k^{2} l
+44\,c_1^{3} c_2^{4} k^{2} l
-48\,c_1^{2} c_2^{5}  k^{2} l
-36\,c_1c_2^{6} k^{2} l
+16\,c_1^{5} c_2k l^{2}
+10\,c_1^{4} c_2^{2} k l^{2}
-16\,c_1^{3} c_2^{3} k l^{2}
+14\,c_1^{2} c_2^{4} k l^{2}
+9\,c_2^{5} l^{2}
+4\,c_1^{7} k^{2}
+c_1^{7} k^{2} l^{2}
+9\,c_2^{7} k^{2} l^{2}
-32\,c_1^{4} c_2^{2} k
-3\,c_1^{4} c_2l^{2}
- 160\,c_1^{3} c_2^{3} k
-54\,c_1^{3} c_2^{2} l^{2}
-96\,c_1^{2} c_2^{4} k
-54\,c_1^{2} c_2^{3} l^{2}
-3\,c_1c_2^{4} l^{2}
-12\,c_1^{6} c_2\,k^{2}
-6\,c_1^{6} k l^{2}
-24\,c_1^{5} c_2^{2} k^{2}
+40\,c_1^{4} c_2^{3} k^{2}
+84\,c_1^{3} c_2^{4} k^{2}
+36\,c_1^{2} c_2^{5} k^{2}
-18\,  c_2^{6} k l^{2}
-48\,c_1^{4} c_2l
-144\,c_1^{3} c_2^{2} l
-144\,c_1^{2} c_2^{3} l
-48\,c_1c_2^{4} l
+4\,c_1^{7} k^{2} l
+32\,c_1^{5} c_2k
-12\,c_1^{6} k l,
\end{autobreak}\\

\begin{autobreak}
NT_{ANB}^k =
(k l^{2}
+4\,k l
+4\,k) c_1^{6}
+(
-6\,c_2k l^{2} 
+4\,c_3k l^{2}
-20\,c_2k l
-4\,c_3k l
-16\,c_2k
-24\,c_3k
-3\,l^{2}
-6\,l) c_1^{5}
+(7\,c_2^{2} k l^{2}
-16\,c_2c_3k l^{2}
+5\,c_3^{2} k l^{2}
+8\,c_2^{2} k l
+32\,c_2c_3k l
-24\,c_3^{2} k l
-8\,c_2^{2} k
+88\,c_2c_3k
+11\,c_2l^{2}
+32\,c_3^{2}  k
-11\,c_3l^{2}
+32\,c_2l
+8\,c_3l
+16\,c_2) c_1^{4}
+(12\,c_2^{3} k l^{2}
-4\,c_2\,c_3^{2} k l^{2}
+56\,c_2^{3} k l
-72\,c_2^{2} c_3k l
+56\,c_2c_3^{2} k l
-8\,c_3^{3} k l
+48\,c_2^{3} k
-8\,c_2^{2} c_3k
-6\,c_2^{2} l^{2}
-160\,c_2c_3^{2} k
+28\,c_2c_3  l^{2}
+24\,c_3^{3} k
-14\,c_3^{2} l^{2}
+4\,c_2^{2} l
+32\,c_2c_3l
+28\,c_3^{2} l
-32\,c_2^{2}
-32\,c_2c_3) c_1^{3}
+(
-17\,c_2^{4} k l^{2}
+40\,c_2^{3} c_3k l^{2}
-42\,c_2^{2} c_3^{2} k l^{2}
+24\,c_2c_3^{3} k l^{2}
-5\,c_3^{4} k l^{2}
-12\,c_2^{4} k l
+56\,c_2^{2} c_3^{2} k l
-64\,c_2c_3^{3} k l
+20\,c_3^{4} k l
+36\,c_2^{4} k
-120\,c_2^{3} c_3k
-2\,c_2^{3} l^{2}
-10\,c_2^{2} c_3l^{2}
+120\,c_2c_3^{3} k
+18\,c_2c_3^{2} l^{2}
-36\,c_3^{4} k
-6\,c_3^{3} l^{2}
+16\,c_2^{3} l
-88\,c_2^{2} c_3l
+32\,c_2c_3^{2} l
+8\,c_3^{3} l
-48\,c_2^{3}
-96\,c_2^{2} c_3
-48\,c_2c_3^{2}) c_1^{2}
+(
-6\,c_2^{5} k l^{2}
-4\,c_2^{4} c_3k l^{2}
+44\,c_2^{3} c_3^{2} k l^{2}
-56\,c_2^{2} c_3^{3} k l^{2}
+26\,c_2c_3^{4} k l^{2}
-4\,c_3^{5} k l^{2}
-36\,c_2^{5} k l
+108\,c_2^{4} c_3k l
-152\,c_2^{3} c_3^{2} k l
+136\,c_2^{2} c_3^{3} k l
-68\,c_2c_3^{4} k l
+12\,c_3^{5} k l
+9\,c_2^{4} l^{2}
+12\,c_2^{3} c_3l^{2}
-2\,c_2^{2} c_3^{2} l^{2}
-4\,c_2c_3^{3} l^{2}
+c_3^{4} l^{2}
+50\,c_2^{4} l
+16\,c_2^{3} c_3l
+4\,c_2^{2} c_3^{2} l
+32\,c_2c_3^{3} l
-6\,c_3^{4} l) c_1
+ 9\,c_2^{6} k l^{2}
-24\,c_2^{5} c_3k l^{2}
+13\,c_2^{4} c_3^{2} k l^{2}
+16\,c_2^{3} c_3^{3} k l^{2}
-21\,c_2^{2} c_3^{4} k l^{2}
+8\,c_2c_3^{5} k l^{2}
-c_3^{6} k l^{2}
-9\,c_2^{5} l^{2}
-3\,c_2^{4} c_3l^{2}
+14\,c_2^{3} c_3^{2} l^{2}
+2\,c_2^{2} c_3^{3} l^{2}
-5\,c_2c_3 ^{4} l^{2}
+c_3^{5} l^{2},
\end{autobreak}\\

\begin{autobreak}
NT_{ANB}^l =
(k^{2} l
+2\,k^{2}) c_1^{7}
+(
-5\,c_2k^{2} l
+3\,c_3k^{2} l
-8\,c_2k^{2}
-4\,c_3k^{2}
-6\,k l
-6\,k) c_1^{6}
+(c_2^{2} k^{2} l
-6\,c_2c_3k^{2} l
+c_3^{2} k^{2} l
-6\,c_2^{2} k^{2}
+24\,c_2c_3k^{2}
-10\,c_3^{2} k^{2}
+16\,c_2k l
-16\,c_3k l
+26\,c_2k
+14\,c_3k
+9\,l) c_1^{5}
+(19\,c_2^{3} k^{2} l
-23\, c_2^{2} c_3k^{2} l
+17\,c_2c_3^{2} k^{2} l
-5\,c_3^{3} k^{2} l
+32\,c_2^{3} k^{2}
-24\,c_2^{2} c_3k^{2}
+8\,c_3^{3} k^{2}
+10\,c_2^{2} k l
+12\,c_2c_3k l
-6\,c_3^{2} k l
+36\,c_2^{2} k
+8\,c_2c_3k
+20\,c_3^{2} k
-3\,c_2l
+21\,c_3l
-24\,c_2) c_1^{4}
+(
-5 \, c_2^{4} k^{2} l
+28\,c_2^{3} c_3k^{2} l
-46\,c_2^{2} c_3^{2} k^{2} l
+28\,c_2c_3^{3} k^{2} l
-5\,c_3^{4} k^{2} l
+22\,c_2^{4} k^{2}
-64\,c_2^{3} c_3k^{2}
+92\,c_2^{2} c_3^{2} k^{2}
-64\,c_2c_3^{3} k^{2}
+14\,c_3^{4} k^{2}
-16\,c_2^{3} k l
+48\,c_2^{2} c_3k l
-48\,  c_2\,c_3^{2} k l
+16\,c_3^{3} k l
+20\,c_2^{3} k
-60\,c_2^{2} c_3k
+28\,c_2c_3^{2} k
-20\,c_3^{3} k
-54\,c_2^{2} l
+4\,c_2c_3l
+10\,c_3^{2} l
-72\,c_2^{2}
-40\,c_2c_3) c_1^{3}
+(
-23\,c_2^{5} k^{2} l
+53\,c_2^{4} c_3k^{2} l
-38\,c_2^{3} c_3^{2} k^{2} l
+10\,  c_2^{2} c_3^{3} k^{2} l
-3\,c_2c_3^{4} k^{2} l
+c_3^{5} k^{2} l
-24\,c_2^{5} k^{2}
+60\,c_2^{4} c_3k^{2}
-48\,c_2^{3} c_3^{2} k^{2}
+8\,c_2^{2} c_3^{3} k^{2}
+8\,c_2c_3^{4} k^{2}
-4\,c_3^{5} k^{2}
+14\,c_2^{4} k l
+8\,c_2^{3} c_3k l
+52\,c_2^{2} c_3^{2} k l
-56\,c_2c_3 ^{3} k l
+14\,c_3^{4} k l
+66\,c_2^{4} k
-88\,c_2^{3} c_3k
+124\,c_2^{2} c_3^{2} k
+8\,c_2c_3^{3} k
-14\,c_3^{4} k
-54\,c_2^{3} l
-34\,c_2^{2} c_3l
+14\,c_2c_3^{2} l
-6\,c_3^{3} l
-72\,c_2^{3}
-80\,c_2^{2} c_3
-8\,c_2c_3^{2}) c_1^{2}
+(3\,c_2^{6} k^{2} l
-22\,c_2 ^{5} c_3k^{2} l
+61\,c_2^{4} c_3^{2} k^{2} l
-84\,c_2^{3} c_3^{3} k^{2} l
+61\,c_2^{2} c_3^{4} k^{2} l
-22\,c_2c_3^{5} k^{2} l
+3\,c_3^{6} k^{2} l
-18\,c_2^{6} k^{2}
+72\,c_2^{5} c_3k^{2}
-130\,c_2^{4} c_3^{2} k^{2}
+144\,c_2^{3} c_3^{3} k^{2}
-102\,c_2^{2} c_3^{4} k^{2} 
+40\,c_2c_3^{5} k^{2}
-6\,c_3^{6} k^{2}
+50\,c_2^{5} k
-34\,c_2^{4} c_3k
-12\,c_2^{3} c_3^{2} k
+28\,c_2^{2} c_3^{3} k
-38\,c_2c_3^{4} k
+6\,c_3^{5} k
-3\,c_2^{4} l
+4\,c_2^{3} c_3l
+14\,c_2^{2} c_3^{2} l
+4\,c_2c_3^{3} l
-3\,c_3^{4} l
-24\,c_2^{4}
-40\,c_2^{3} c_3
-8\,  c_2^{2} c_3^{2}
+8\,c_2c_3^{3}) c_1
+c_3^{5} l
+9\,c_2^{5} l
-33\,c_2^{6} c_3k^{2} l
-9\,c_2c_3^{6} k^{2} l
+12\,c_2^{5} c_3k l
+34\,c_2^{4} c_3^{2} k l
-24\,c_2^{3} c_3^{3} k l
-14\,c_2^{2} c_3^{4} k l
+12\,c_2c_3^{5} k l
+37\,c_2^{5} c_3^{2} k^{2} l
+3\,c_2^{4} c_3^{3} k^{2} l
-37\,c_2^{3} c_3^{4} k^{2} l
+29\,c_2^{2} c_3^{5} k^{2} l
+c_3^{7} k^{2} l
+21\,c_2^{4} c_3\,l
-6\,c_2^{2} c_3^{3} l
-3\,c_2c_3^{4} l
-18\,c_2^{6} k l
-2\,c_3^{6} k l
+9\,c_2^{7} k^{2} l
+10\,c_2^{3} c_3^{2} l,
\end{autobreak}\\

\begin{autobreak}
NT_{ANB}^{c_1} =
(k^{2} l
+2\,k^{2}) c_1^{7}
+(
-3\,c_2k^{2} l
+l k ^{2} c_3
-4\,c_2k^{2}
-8\,k^{2} c_3
-4\,k l
-2\,k) c_1^{6}
+(
-7\,c_2^{2} k^{2} l
+6\,c_2c_3k^{2} l
-3\,l k^{2} c_3^{2}
-18\,c_2^{2} k^{2}
+24\,c_2c_3k^{2}
+2\,k^{2} c_3^{2}
+6\,c_2k l
-6\,l k c_3
+22\,c_2k
+2\,k c_3
+3\,l) c_1^{5}
+(13\,c_2^{3} k^{2} l
-17 \, c_2^{2} c_3k^{2} l
+15\,c_2\,c_3^{2} k^{2} l
-3\,c_3^{3} k^{2} l
+8\,c_2^{3} k^{2}
+48\,c_2^{2} c_3k^{2}
-56\,c_2c_3^{2} k^{2}
+16\,c_3^{3} k^{2}
+18\,c_2^{2} k l
-4\,c_2c_3k l
+2\,c_3^{2} k l
+28\,c_2^{2} k
-8\,c_2c_3k
+12\,c_3^{2} k
-5\,c_2l
+11\,c_3l
-8\,  c_2) c_1^{4}
+(27\,c_2^{4} k^{2} l
-44\,c_2^{3} c_3k^{2} l
+26\,c_2^{2} c_3^{2} k^{2} l
-12\,c_2c_3^{3} k^{2} l
+3\,c_3^{4} k^{2} l
+62\,c_2^{4} k^{2}
-80\,c_2^{3} c_3k^{2}
+12\,c_2^{2} c_3^{2} k^{2}
+16\,c_2\,c_3^{3} k^{2}
-10\,c_3^{4} k^{2}
-4\,c_2^{3} k l
+12\,c_2^{2}  c_3k l
-12\,c_2c_3^{2} k l
+4\,c_3^{3} k l
-68\,c_2^{3} k
+12\,c_2^{2} c_3k
-44\,c_2c_3^{2} k
+4\,c_3^{3} k
-10\,c_2^{2} l
-12\,c_2c_3l
+14\,c_3^{2} l
+8\,c_2^{2}
-24\,c_2c_3) c_1^{3}
+(
-c_2^{5} k^{2} l
+15\,c_2^{4} c_3k^{2} l
-42\,c_2^{3} c_3^{2} k^{2}  l
+46\,c_2^{2} c_3^{3} k^{2} l
-21\,c_2c_3^{4} k^{2} l
+3\,c_3^{5} k^{2} l
+60\,c_2^{5} k^{2}
-168\,c_2^{4} c_3k^{2}
+200\,c_2^{3} c_3^{2} k^{2}
-144\,c_2^{2} c_3^{3} k^{2}
+60\,c_2c_3^{4} k^{2}
-8\,c_3^{5} k^{2}
-122\,c_2^{4} k
+56\,c_2^{3} c_3k
-44\,c_2^{2} c_3^{2} k 
+24\,c_2c_3^{3} k
-10\,c_3^{4} k
-2\,c_2^{3} l
-14\,c_2^{2} c_3l
-6\,c_2c_3^{2} l
+6\,c_3^{3} l
+40\,c_2^{3}
+16\,c_2^{2} c_3
-24\,c_2c_3^{2}) c_1^{2}
+(
-21\,c_2^{6} k^{2} l
+70\,c_2^{5} c_3k^{2} l
-87\,c_2^{4} c_3^{2} k^{2} l
+52\,c_2^{3} c_3^{3} k^{2} l
-19\,c_2 ^{2} c_3^{4} k^{2} l
+6\,c_2c_3^{5} k^{2} l
-c_3^{6} k^{2} l
+18\,c_2^{6} k^{2}
-72\,c_2^{5} c_3k^{2}
+130\,c_2^{4} c_3^{2} k^{2}
-144\,c_2^{3} c_3^{3} k^{2}
+102\,c_2^{2} c_3^{4} k^{2}
-40\,c_2c_3^{5} k^{2}
+6\,c_3^{6} k^{2}
+30\,c_2^{5} k l
-22\,c_2^{4} c_3k l
-20\,  c_2^{3} c_3^{2} k l
+20\,c_2^{2} c_3^{3} k l
-10\,c_2c_3^{4} k l
+2\,c_3^{5} k l
-50\,c_2^{5} k
+34\,c_2^{4} c_3k
+12\,c_2^{3} c_3^{2} k
-28\,c_2^{2} c_3^{3} k
+38\,c_2c_3^{4} k
-6\,c_3^{5} k
-9\,c_2^{4} l
-12\,c_2^{3} c_3l
+2\,c_2^{2} c_3^{2} l
+4\,c_2c_3^{3} l 
-c_3^{4} l
+24\,c_2^{4}
+40\,c_2^{3} c_3
+8\,c_2^{2} c_3^{2}
-8\,c_2c_3^{3}) c_1
-9\,c_2^{5} l
-c_3^{5} l
+33\,c_2^{6} c_3k^{2} l
-37\,c_2^{5} c_3^{2} k^{2} l
-3\,c_2^{4} c_3^{3} k^{2} l
+37\,c_2^{3} c_3^{4} k^{2} l
-29\,c_2^{2} c_3^{5} k^{2} l
+9\,c_2c_3^{6} k^{2} l
-12\,  c_2^{5} c_3k l
-34\,c_2^{4} c_3^{2} k l
+24\,c_2^{3} c_3^{3} k l
+14\,c_2^{2} c_3^{4} k l
-12\,c_2\,c_3^{5} k l
-9\,c_2^{7} k^{2} l
-c_3^{7} k^{2} l
+18\,c_2^{6} k l
+2\,c_3^{6} k l
-21\,c_2^{4} c_3l
-10\,c_2^{3} c_3^{2} l
+6\,c_2^{2} c_3^{3} l
+3\,c_2c_3^{4} l,
\end{autobreak}\\

\begin{autobreak}
NT_{ANB}^{c_2} =
(k^{2} l^{2}
+4\,k^{2} l
+4\,k^{2}) c_1^{7}
+(
-c_2k ^{2} l^{2}
+3\,c_3k^{2} l^{2}
-8\,c_3k^{2} l
+4\,c_2k^{2}
-28\,c_3k^{2}
-6\,k l^{2}
-12\,k l) c_1^{6}
+(
-7\,c_2^{2} k^{2} l^{2}
+2\,c_2c_3k^{2} l^{2}
+c_3^{2} k^{2} l^{2}
-28\,c_2^{2} k^{2} l
+16\,c_2c_3k^{2} l
-20\,c_3^{2} k^{2} l
-24\,c_2^{2} k^{2}
-16\,c_2  c_3k^{2}
+56\,c_3^{2} k^{2}
-16\,c_3k l^{2}
+4\,c_2k l
+28\,c_3k l
+16\,c_2k
+9\,l^{2}) c_1^{5}
+(
-c_2^{3} k^{2} l^{2}
-15\,c_2^{2} c_3k^{2} l^{2}
+13\,c_2c_3^{2} k^{2} l^{2}
-5\,c_3^{3} k^{2} l^{2}
-32\,c_2^{3} k^{2} l
+48\,c_2^{2} c_3k^{2} l
-32\,c_2c_3^{2}  k^{2} l
+16\,c_3^{3} k^{2} l
-56\,c_2^{3} k^{2}
+120\,c_2^{2} c_3k^{2}
+18\,c_2^{2} k l^{2}
-24\,c_2c_3^{2} k^{2}
-12\,c_2c_3k l^{2}
-8\,c_3^{3} k^{2}
-6\,c_3^{2} k l^{2}
+56\,c_2^{2} k l
+48\,c_2c_3k l
+40\,c_3^{2} k l
+48\,c_2^{2} k
-48\,c_2c_3k
+9\,c_2  l^{2}
+21\,c_3l^{2}
-24\,c_2l) c_1^{4}
+(11\,c_2^{4} k^{2} l^{2}
-20\,c_2^{3} c_3 k^{2} l^{2}
+2\,c_2^{2} c_3^{2} k^{2} l^{2}
+12\,c_2c_3^{3} k^{2} l^{2}
-5\,c_3^{4} k^{2} l^{2}
+12\,c_2^{4} k^{2} l
-32\,c_2^{3} c_3k^{2} l
+56\,c_2^{2} c_3^{2} k^{2} l
-64\,c_2c_3^{3}  k^{2} l
+28\,c_3^{4} k^{2} l
-44\,c_2^{4} k^{2}
+176\,c_2^{3} c_3k^{2}
-216\,c_2^{2} c_3^{2} k^{2}
+16\,c_2^{2} c_3k l^{2}
+144\,c_2c_3^{3} k^{2}
-32\,c_2c_3^{2} k l^{2}
-60\,c_3^{4} k^{2}
+16\,c_3^{3} k l^{2}
+24\,c_2^{3} k l
+24\,c_2^{2} c_3k l
+56\,c_2c_3^{2} k l
-40\,c_3^{3} k l
+48\,c_2^{3} k
-96\,c_2^{2} c_3k
+2\,c_2^{2} l^{2}
-16\,c_2c_3^{2} k
+12\,c_2c_3l^{2}
+10\,c_3^{2} l^{2}
-40\,c_2^{2} l
-40\,c_2c_3l) c_1^{3}
+(5\,c_2^{5} k^{2} l^{2}
+5\,c_2^{4} c_3k^{2} l^{2}
-30\,c_2^{3} c_3^{2} k^{2} l^{2}
+26\,c_2^{2}  c_3^{3} k^{2} l^{2}
-7\,c_2c_3^{4} k^{2} l^{2}
+c_3^{5} k^{2} l^{2}
+32\,c_2^{5} k^{2} l
-104\,c_2^{4} c_3k^{2} l
+128\,c_2^{3} c_3^{2} k^{2} l
-80\,c_2^{2} c_3^{3} k^{2} l
+32\,c_2c_3^{4} k^{2} l
-8\,c_3^{5} k^{2} l
-12\,c_2^{5} k^{2}
+68\,c_2^{4} c_3k^{2}
-18\,c_2^{4} k l^{2}
-136\,c_2^{3} c_3^{2} k^{2}
+8\,c_2^{3} c_3k l^{2}
+152\,c_2^{2} c_3^{3} k^{2}
-12\,c_2^{2} c_3^{2} k l^{2}
-108\,c_2c_3^{4} k^{2}
-24\,c_2c_3^{3} k l^{2}
+36\,c_3^{5} k^{2}
+14\,c_3^{4} k l^{2}
-44\,c_2^{4} k l
+16\,c_2^{3} c_3k l
+72\,c_2^{2} c_3^{2} k  l
-16\,c_2c_3^{3} k l
-28\,c_3^{4} k l
+16\,c_2^{4} k
-48\,c_2^{3} c_3k
+10\,c_2^{3} l^{2}
-16\,c_2^{2} c_3^{2} k
+14\,c_2^{2} c_3l^{2}
+48\,c_2c_3^{3} k
-2\,c_2c_3^{2} l^{2}
-6\,c_3^{3} l^{2}
-8\,c_2^{3} l
-16\,c_2^{2} c_3l
-8\,c_2c_3^{2} l) c_1^{2}
+(
-5\,c_2 ^{6} k^{2} l^{2}
+18\,c_2^{5} c_3k^{2} l^{2}
-19\,c_2^{4} c_3^{2} k^{2} l^{2}
-4\,c_2^{3} c_3^{3} k^{2} l^{2}
+21\,c_2^{2} c_3^{4} k^{2} l^{2}
-14\,c_2c_3^{5} k^{2} l^{2}
+3\,c_3^{6} k^{2} l^{2}
+12\,c_2^{6} k^{2} l
-48\,c_2^{5} c_3k^{2} l
+60\,c_2^{4} c_3^{2} k^{2} l
-60  c_2^{2} c_3^{4} k^{2} l
+48\,c_2c_3^{5} k^{2} l
-12\,c_3^{6} k^{2} l
-28\,c_2^{5} k l
+12\,c_2^{4} c_3k l
+56\,c_2^{3} c_3^{2} k l
-24\,c_2^{2} c_3^{3} k l
-28\,c_2c_3^{4} k l
+12\,c_3^{5} k l
+5\,c_2^{4} l^{2}
+12\,c_2^{3} c_3l^{2}
+6\,c_2^{2} c_3^{2} l^{2}
-4\,  c_2c_3^{3} l^{2}
-3\,c_3^{4} l^{2}
+8\,c_2^{4} l
+24\,c_2^{3} c_3l
+24\,c_2^{2} c_3^{2} l
+8\,c_2\,c_3^{3} l) c_1
-11\,c_2^{4} c_3l^{2}
-14\,c_2^{3} c_3^{2} l^{2}
-6\,c_2^{2} c_3^{3} l^{2}
+c_2c_3^{4} l^{2}
+6\,c_2^{6} k l^{2}
-3\,c_2^{7} k^{2} l^{2}
-3\,c_2^{5} l^{2}
+4\,c_2 ^{5} c_3k l^{2}
-14\,c_2^{4} c_3^{2} k l^{2}
-8\,c_2^{3} c_3^{3} k l^{2}
+10\,c_2^{2} c_3^{4} k l^{2}
+4\,c_2c_3^{5} k l^{2}
+c_3^{5} l^{2}
-2\,c_3^{6} k l^{2}
+c_3^{7} k^{2} l^{2}
+7\,c_2^{6} c_3k^{2} l^{2}
+c_2^{5} c_3^{2} k^{2} l^{2}
-13\,c_2^{4} c_3^{3} k^{2} l^{2}
+7\,  c_2^{3} c_3^{4} k^{2} l^{2}
+5\,c_2^{2} c_3^{5} k^{2} l^{2}
-5\,c_2c_3^{6} k^{2} l^{2},
\end{autobreak}\\

\begin{autobreak}
NT_{ANB}^{c_3} =
(3\,k l^{2}
+6\,k l) c_1^{6}
+(
-16\,l^{2} k c_2
+ 8\,c_3k l^{2}
-34\,c_2k l
-14\,c_3k l
-16\,c_2k
-9\,l^{2}) c_1^{5}
+(11\,c_2^{2} k l^{2}
-22\,l^{2} k c_3c_2
+3\,c_3^{2} k l^{2}
+12\,c_2^{2} k l
+56\,c_2k l c_3
-20\,c_3^{2} k l
+16\,c_2^{2} k
+48\,c_2k c_3
+27\,l^{2} c_2
-21\,c_3l^{2}
+72\,c_2l ) c_1^{4}
+(40\,c_2^{3} k l^{2}
-24\,c_2^{2} c_3k l^{2}
+24\,c_2c_3^{2} k l^{2}
-8\,c_3^{3} k l^{2}
+108\,c_2^{3} k l
-52\,c_2^{2} c_3k l
-12\,c_2c_3^{2} k l
+20\,c_3^{3} k l
+80\,c_2^{3} k
-160\,c_2^{2} c_3k
-2\,c_2^{2} l^{2}
+16\,c_2c_3^{2} k
+36\,c_2  c_3l^{2}
-10\,c_3^{2} l^{2}
-40\,c_2^{2} l
+120\,c_2c_3l
-128\,c_2^{2}) c_1^{3}
+(
-23\,c_2^{4} k l^{2}
+44\,c_2^{3} c_3k l^{2}
-74\,c_2^{2} c_3^{2} k l^{2}
+44\,c_2c_3^{3} k l^{2}
-7\,c_3^{4} k l^{2}
+14\,c_2^{4} k l
-72\,c_2^{3} c_3k l
+84\,c_2^{2} c_3^{2} k  l
-72\,c_2c_3^{3} k l
+14\,c_3^{4} k l
+48\,c_2^{4} k
-208\,c_2^{3} c_3k
-2\,c_2^{3} l^{2}
+208\,c_2^{2} c_3^{2} k
-14\,c_2^{2} c_3l^{2}
-48\,c_2c_3^{3} k
-6\,c_2c_3^{2} l^{2}
+6\,c_3^{3} l^{2}
-40\,c_2^{3} l
-16\,c_2^{2} c_3l
+24\,c_2c_3^{2} l
-128\,c_2^{3}
- 128\,c_2^{2} c_3) c_1^{2}
+(
-24\,c_2^{5} k l^{2}
+32\,c_2^{4} c_3k l^{2}
+16\,c_2^{3} c_3^{2} k l^{2}
-32\,c_2^{2} c_3^{3} k l^{2}
+8\,c_2c_3^{4} k l^{2}
-42\,c_2^{5} k l
+50\,c_2^{4} c_3k l
+12\,c_2^{3} c_3^{2} k l
-44\,c_2^{2} c_3^{3} k l
+30\,c_2c_3^{4} k  l
-6\,c_3^{5} k l
+27\,c_2^{4} l^{2}
+36\,c_2^{3} c_3l^{2}
-6\,c_2^{2} c_3^{2} l^{2}
-12\,c_2c_3^{3} l^{2}
+3\,c_3^{4} l^{2}
+72\,c_2^{4} l
+120\,c_2^{3} c_3l
+24\,c_2^{2} c_3^{2} l
-24\,c_2c_3^{3} l) c_1
+9\,c_2^{6} k l^{2}
-6\,c_2^{5} c_3k l^{2}
-17\,c_2^{4} c_3^{2}  k l^{2}
+12\,c_2^{3} c_3^{3} k l^{2}
+7\,c_2^{2} c_3^{4} k l^{2}
-6\,c_2c_3^{5} k l^{2}
+c_3^{6} k l^{2}
-9\,c_2^{5} l^{2}
-21\,c_2^{4} c_3l^{2}
-10\,c_2^{3} c_3^{2} l^{2}
+6\,c_2^{2} c_3^{3} l^{2}
+3\,c_2c_3^{4} l^{2}
-c_3^{5} l^{2},
\end{autobreak}\\

\begin{autobreak}
NS_{LNB} =
4\,k^{2} c_3^{6}
+(
-8\,c_1k^{2}
+4\,k) c_3^{5}
+(
-4\,c_1^{2} k^{2}
+16\,c_1c_2k^{2}
-72\,c_2^{2} k^{2}
-12\,c_1k
+1) c_3^{4}
+(16\,c_1^{3} k^{2}
-32\,c_1^{2} c_2k^{2}
+112\,c_1c_2^{2} k^{2}
+128\,c_2^{3} k^{2}
-8\,c_1^{2} k
+8\,c_1c_2k
+48\,c_2^{2} k
-4\,c_1
+4\,c_2) c_3^{3}
+(
-4\,c_1^{4} k^{2}
-16\,c_1^{2} c_2^{2} k^{2}
-144\,c_1  c_2^{3} k^{2}
-60\,c_2^{4} k^{2}
+24\,c_1^{3} k
-56\,c_1^{2} c_2k
+136\,c_1c_2^{2} k
-8\,c_2^{3} k
-2\,c_1^{2}
+4\,c_1c_2
-10\,c_2^{2}) c_3^{2}
+(
-8\,c_1^{5} k^{2}
+32\,c_1^{4} c_2k^{2}
-16\,c_1^{3} c_2^{2} k^{2}
-32\,c_1^{2} c_2^{3} k^{2}
+24\,c_1c_2^{4} k^{2}
+4\,c_1^{4} k
-8\,c_1^{3} c_2k
-80\,c_1^{2} c_2^{2} k
-248\,c_1c_2^{3} k
-52\,c_2^{4} k
+12\,c_1^{3}
-36\,c_1^{2} c_2
-76\,c_2^{2} c_1
-28\,c_2^{3}) c_3
+4\,c_1^{6} k^{2}
-16\,c_1^{5} c_2k^{2}
-8\,c_1^{4} c_2^{2} k^{2}
+48\,c_1^{3} c_2^{3} k^{2}
+36\,c_1^{2} c_2^{4} k^{2}
-12\,c_1^{5} k
+56\,c_1^{4} c_2k
+24\,c_1^{3} c_2^{2} k
-32\,c_1^{2} c_2^{3} k
+20\,c_1c_2^{4} k
+8\,c_2^{5} k
+9\,c_1^{4}
-36\,c_1^{3} c_2
-114\,c_1^{2} c_2^{2}
-84\,c_1c_2^{3}
-15\,c_2^{4},
\end{autobreak}\\

\begin{autobreak}
NT_{LNB}^k =
2\,k c_1^{5}
+(
-10\,c_2k
-2\,c_3k
-3) 
c_1^{4}
+(6\,c_2^{2} k
+8\, c_2c_3k
-4\,c_3^{2} k
+17\,c_2
-2\,c_3) c_1^{3}
+(18\,c_2^{3} k
-10\,c_2^{2} c_3k
+12\,c_2c_3^{2} k
+4\,k c_3^{3}
-11\,c_2^{2}
+17\,c_3c_2
+4\,c_3^{2}) c_1^{2}
+(12\,c_2^{3} c_3k
-30\,c_2^{2} c_3^{2} k
-8\,c_2c_3^{3} k
+2\,c_3^{4} k
+3\,c_2^{3}
-48\,c_2^{2} c_3
- c_2c_3^{2}
+2\,c_3^{3}) c_1
-30\,c_2^{3} c_3^{2} k
+34\,c_2^{2} c_3^{3} k
-2\,c_2c_3^{4} k
-2\,c_3^{5} k
+2\,c_2^{4}
-11\,c_2^{3} c_3
-13\,c_2^{2} c_3^{2}
-c_2\,c_3^{3}
-c_3^{4},
\end{autobreak}\\

\begin{autobreak}
NT_{LNB}^{c_1} =
2\,k^{2} c_1^{5}
+(
-8\,k^{2} c_2
-2\,k^{2} c_3
-7\,k) c_1^{4}
+(
-4\,c_2 ^{2} k^{2}
+16\,k^{2} c_3\,c_2
-4\,k^{2} c_3^{2}
+28\,k c_2
-2\,c_3k
+6) c_1^{3}
+(24\,c_2^{3} k^{2}
-20\,c_2^{2} c_3k^{2}
+4\,c_3^{3} k^{2}
+14\,c_2^{2} k
-6\,k c_3c_2
+12\,c_3^{2} k
-18\,c_2
+10\,c_3) c_1^{2}
+(18\,c_2^{4} k^{2}
-32\,c_2^{3} c_3k^{2}
+28\,c_2^{2} c_3 ^{2} k^{2}
-16\,c_2c_3^{3} k^{2}
+2\,c_3^{4} k^{2}
-20\,c_2^{3} k
+26\,c_2^{2} c_3k
-32\,c_2c_3^{2} k
+2\,k c_3^{3}
-38\,c_2^{2}
-20\,c_3c_2
+2\,c_3^{2}) c_1
+6\,c_2^{4} c_3k^{2}
-8\,c_2^{3} c_3^{2} k^{2}
-4\,c_2^{2} c_3^{3} k^{2}
+8\,c_2c_3^{4} k^{2}
-2\,c_3^{5} k^{2}
+c_2 ^{4} k
-34\,c_2^{3} c_3k
+36\,c_2^{2} c_3^{2} k
+2\,c_2c_3^{3} k
-5\,c_3^{4} k
-14\,c_2^{3}
-14\,c_2^{2} c_3
-2\,c_2c_3^{2}
-2\,c_3^{3},
\end{autobreak}\\

\begin{autobreak}
NT_{LNB}^{c_2} =
4\,k^{2} c_1^{5}
+(
-4\,c_2k^{2}
-12\,c_3k^{2}
-12\,k) c_1^{4}
+(
-20  c_2^{2} k^{2}
+16\,c_2c_3k^{2}
+8\,c_3^{2} k^{2}
+24\,k c_2
+16\,c_3k
+9) c_1^{3}
+(
-12\,c_2^{3} k^{2}
+36\,c_2^{2} c_3k^{2}
-24\,c_2c_3^{2} k^{2}
+8\,c_3^{3} k^{2}
+32\,c_2^{2} k
-32\,c_2c_3k
+8\,c_3^{2} k
-21\,c_2
+3\,c_3) c_1^{2}
+(8\,c_2^{3} c_3k^{2}
-12  c_2^{2} c_3^{2} k^{2}
+16\,c_2c_3^{3} k^{2}
-12\,c_3^{4} k^{2}
-8\,c_2^{3} k
+8\,c_2c_3^{2} k
-16\,k c_3^{3}
+3\,c_2^{2}
-2\,c_3c_2
-5\,c_3^{2}) c_1
+4\,c_2^{3} c_3^{2} k^{2}
-4\,c_2^{2} c_3^{3} k^{2}
-4\,c_2c_3^{4} k^{2}
+4\,c_3^{5} k^{2}
-4\,c_2^{4} k
+4\,c_3^{4} k
+c_2^{3}
+3 \, c_2^{2} c_3
+3\,c_2c_3^{2}
+c_3^{3},
\end{autobreak}\\

\begin{autobreak}
NT_{LNB}^{c_3} =
6\,k c_1^{5}
+(
-40\,c_2k
-2\,k c_3
-9) c_1^{4}
+(36\,c_2^{2} k
+28\,  c_2c_3k
-12\,c_3^{2} k
+72\,c_2
-12\,c_3) c_1^{3}
+(88\,c_2^{3} k
-96\,c_2^{2} c_3k
+44\,c_2c_3^{2} k
+4\,k c_3^{3}
-90\,c_2^{2}
+72\,c_3c_2
+2\,c_3^{2}) c_1^{2}
+(6\,c_2^{4} k
-28\,c_2^{3} c_3k
-4\,c_2^{2} c_3^{2} k
-28\,c_2\,c_3^{3} k
+6\,c_3^{4} k
-48\,c_2^{3} 
-60\,c_2^{2} c_3
-8\,c_2\,c_3^{2}
+4\,c_3^{3}) c_1
-30\,c_2^{4} c_3k
+4\,c_2^{3} c_3^{2} k
+32\,c_2^{2} c_3^{3} k
-4\,c_2c_3^{4} k
-2\,c_3^{5} k
-5\,c_2^{4}
-16\,c_2^{3} c_3
-18\,c_2^{2} c_3^{2}
-8\,c_2c_3^{3}
-c_3^{4}.
\end{autobreak}\\

\end{align*}

\end{CJK}
\end{document}